\documentclass[10pt]{amsart}

\usepackage{fullpage}
\usepackage{amssymb}
\usepackage{amsmath}
\usepackage{amsxtra}
\usepackage{amscd}
\usepackage{graphicx}
\usepackage{amsfonts}
\usepackage{pb-diagram}
\usepackage{float}
\usepackage{enumerate}
\usepackage{dsfont}

\numberwithin{equation}{section}
\newtheorem{thm}[equation]{Theorem}

\newtheorem{cor}[equation]{Corollary}
\newtheorem{prop}[equation]{Proposition}

\newtheorem{defn}[equation]{Definition}
\newtheorem{conj}[equation]{Conjecture}
\newtheorem{exmp}[equation]{Example}

\newtheorem{rem}[equation]{Remark}

\newtheorem{hyp}[equation]{Hypothesis}

\def\leq{\leqslant}

\def\Q{{\mathbb Q}}

\def\Z{{\mathbb Z}}
\def\R{{\mathbb R}}
\def\C{{\mathbb C}}
\def\K{{\mathbb K}}
\def\Irr{{\mathrm{Irr}}}
\def\KZ{{\mathrm{KZ}}}
\def\bc{{\boldsymbol{\rm c}}}

\def\bl{{\boldsymbol{\lambda}}}

\begin{document}

\title[Hecke algebras and symplectic reflection algebras]
  {Hecke algebras and symplectic reflection algebras}
\author{Maria Chlouveraki}
\address{Laboratoire de Math\'ematiques UVSQ, B\^atiment Fermat, 45 avenue des \'Etats-Unis,  78035 Versailles cedex, France.}
\email{maria.chlouveraki@uvsq.fr}

\maketitle

\maketitle

\begin{abstract}
The current article is a short survey on the theory of Hecke algebras, and in particular Kazhdan--Lusztig theory, and on the theory 
of symplectic reflection algebras, and in particular rational Cherednik algebras. The emphasis is on the connections
between Hecke algebras and rational Cherednik algebras that could allow us to obtain a generalised Kazhdan--Lusztig theory, or at least its applications, for all complex reflection groups.
\end{abstract}

\tableofcontents

\section{Introduction}

Finite Coxeter groups are finite groups of real matrices that are generated by reflections. 
 They include the Weyl groups, which are  fundamental  in the classification of simple complex Lie algebras as well as simple algebraic groups. 
Iwahori--Hecke algebras associated to Weyl groups appear naturally as endomorphism algebras of  induced representations in the study of finite reductive groups.
They can also be defined independently as deformations of group algebras of finite Coxeter groups, where the deformation depends on an indeterminate $q$ and a weight function $L$. For $q=1$, we recover the group algebra. For a finite Coxeter group $W$, we will denote by $\mathcal{H}(W,L)$ the associated Iwahori--Hecke algebra. 

When $q$ is an indeterminate, the Iwahori--Hecke algebra $\mathcal{H}(W,L)$ is semisimple.
By Tits's deformation theorem, there exists a bijection between the set of irreducible representations of $\mathcal{H}(W,L)$ and the set  $\Irr(W)$ of  irreducible representations of $W$. Using this bijection, Lusztig attaches to every irreducible representation of $W$ an integer depending on $L$, thus defining the famous $a$\emph{-function}. The $a$-function is used in his definition of \emph{families of characters}, a partition of $\Irr(W)$ which plays a key
role in the organisation of families of unipotent characters in the case of finite reductive groups. 

Kazhdan--Lusztig theory is a key to understanding the representation theory of the Iwahori--Hecke algebra $\mathcal{H}(W,L)$. There exists a special basis of $\mathcal{H}(W,L)$, called the \emph{Kazhdan--Lusztig basis}, which allows us to define the \emph{Kazhdan--Lusztig cells} for  $\mathcal{H}(W,L)$, a certain set of equivalence classes on $W$.
The construction of Kazhdan--Lusztig cells yields the construction of representations for $\mathcal{H}(W,L)$. It also gives another, more combinatorial, definition for Lusztig's families of characters.

Now, when $q$ specialises to a non-zero complex number $\eta$, and more specifically to a root of unity, the specialised Iwahori--Hecke algebra $\mathcal{H}_\eta(W,L)$  is not necessarily semisimple and we no longer have a bijection between its irreducible representations and $\Irr(W)$. We obtain then a decomposition matrix which records  how the irreducible representations of the semisimple algebra split after the specialisation. A \emph{canonical basic set} is a subset of $\Irr(W)$ in bijection with the irreducible representations of $\mathcal{H}_\eta(W,L)$ (and thus a labelling set for the columns of the decomposition matrix) with good properties. Its good properties ensure that the decomposition matrix has a lower unitriangular form while the $a$-function increases (roughly) down the columns. Canonical basic sets were defined by Geck and Rouquier \cite{GR}, who also proved their existence in certain cases with the use of Kazhdan--Lusztig theory. Thanks to the work of many people, canonical basic sets are now proved to exist and explicitly described for all finite Coxeter groups and for any choice of $L$.

Finite Coxeter groups are particular cases of complex reflection groups, that is, finite groups of complex matrices  generated by  ``pseudo-reflections''. Their classification is due to Shephard and Todd \cite{ShTo}: An irreducible complex reflection group either belongs to the infinite series $G(\ell,p,n)$ or is one of the 34 exceptional groups $G_4,\ldots,G_{37}$  (see Theorem \ref{ShToClas}). Important work in the last two decades has suggested that complex reflection groups will play a crucial, but not yet understood role in representation theory, and may even become as ubiquitous in the study of other mathematical structures. In fact, they behave so much  like real reflection groups that Brou\'e, Malle and Michel \cite{BMM2} conjectured that they could play the role of Weyl groups for, as yet mysterious, objects generalising finite reductive groups. These objects are called ``Spetses''.

Brou\'e, Malle and Rouquier \cite{BMR} defined Hecke algebras for complex reflection groups as deformations of their group algebras.
A generalised Kazhdan--Lusztig cell theory for these algebras, known as \emph{cyclotomic Hecke algebras}, is expected to help find Spetses. Unfortunately, we do not have a Kazhdan--Lusztig basis for complex reflection groups. However, we can define families of characters using Rouquier's definition:
In \cite{Ro} Rouquier gave an alternative definition for Lusztig's families of characters by proving that, in the case of Weyl groups, they coincide with the blocks of the Iwahori--Hecke algebra over a certain ring, called the \emph{Rouquier ring}. This definition generalises without problem to the case of complex reflection groups and their cyclotomic Hecke algebras, producing the so-called \emph{Rouquier families}. These families have now been determined for all cyclotomic Hecke algebras of all complex reflection groups, see \cite{mybook}.

We also have an $a$-function and can define canonical basic sets for cyclotomic Hecke algebras. Although there is no Kazhdan--Lusztig theory in the complex case, canonical basic sets are now known to exist for the groups of the infinite series $G(\ell,p,n)$ and for some exceptional ones. In order to obtain canonical basic sets for $G(\ell,1,n)$, Geck and Jacon used Ariki's Theorem on the categorification of Hecke algebra representations and Uglov's work on canonical bases for higher level Fock spaces  \cite{GJ, jac1, jaca, GJlivre}. The result  for $G(\ell,p,n)$ derives from that for $G(\ell,1,n)$ with the use of Clifford Theory  \cite{geja, chja2}. 

In this paper we will  see how we could use the representation theory of symplectic reflection algebras, and in particular rational Cherednik algebras, to obtain families of characters and canonical basic sets for cyclotomic Hecke algebras associated with complex reflection groups.

Symplectic reflection algebras are related to a large number of areas of mathematics such as combinatorics, integrable systems, real algebraic geometry, quiver varieties, symplectic  resolutions of singularities and, of course, representation theory. 
They were introduced by Etingof and Ginzburg in \cite{EG} for the study of symplectic resolutions of the orbit space $V/G$, where $V$ is a symplectic complex vector space and $G \subset {\rm Sp(V)}$ is a finite group acting on $V$. Verbitsky \cite{ver} has shown that $V/G$ admits a symplectic resolution only if $(G,V)$ is a symplectic reflection group, that is, $G$ is generated by symplectic reflections. Thanks to the insight by Etingof and Ginzburg, the study of the representation theory of symplectic reflection algebras has led to the (almost) complete classification of symplectic reflection groups $(G,V)$ such that $V/G$ admits a symplectic resolution.

Let $(G,V)$ be a  symplectic reflection group, and  let $TV^{*}$ denote the tensor algebra on the dual space $V^*$ of $V$. The symplectic reflection algebra ${\bf H}_{t,\bc}(G)$ associated to $(G,V)$ is defined as the quotient of  $TV^{*} \rtimes G$ by certain relations depending on a complex function $\bc$ and a parameter $t$. The representation theory of ${\bf H}_{t,\bc}(G)$  varies a lot according to whether  $t$ is zero or not.  A complex reflection group $W \subset {\rm GL}(\mathfrak{h})$, where $\mathfrak{h}$ is a complex vector space, can be seen as a symplectic reflection group acting on $V =  \mathfrak{h}\, \oplus\, \mathfrak{h}^*$. Symplectic reflection algebras associated with complex reflection groups are known as \emph{rational Cherednik algebras}. 

If $t \neq 0$, there exists an important category of representations of the rational Cherednik algebra, the \emph{category} $\mathcal{O}$, and an exact functor, the $\KZ$-\emph{functor}, from $\mathcal{O}$ to the category of representations of a certain specialised cyclotomic Hecke algebra $\mathcal{H}_\eta(W)$ (the specialisation depends on the choice of parameters for the rational Cherednik algebra --- every specialised Hecke algebra can arise this way). Category $\mathcal{O}$ is a highest weight category, and it comes equipped with a set of standard modules $\{\Delta(E)\,|\, E \in \Irr(W)\}$, a set of simple modules $\{{\rm L}(E) \,|\, E \in \Irr(W)\}$  and a decomposition matrix that records the number of times that ${\rm L}(E)$ appears in the composition series of $\Delta(E')$ for $E,E' \in \Irr(W)$. The exactness of $\KZ$ allows us to read off the decomposition matrix of $\mathcal{H}_\eta(W)$ from the decomposition matrix of category $\mathcal{O}$. Using this, we proved in \cite{CGG} the existence of canonical basic sets for all finite Coxeter groups and for complex reflection groups of type $G(\ell,1,n)$. In particular, we showed that $E$ belongs to the canonical basic set for $\mathcal{H}_\eta(W)$ if and only if 
$\KZ({\rm L}(E)) \neq 0$. Our proof of existence is quite general and it does not make use of Ariki's Theorem for type $G(\ell,1,n)$. However, the explicit description of canonical basic sets in these cases by previous works answers simultaneously the question of which simple modules are killed by the $\KZ$-functor; this appears to be new. We also proved that the images of the standard modules via the $\KZ$-functor  are isomorphic to the cell modules of Hecke algebras with cellular structure, but we will not go into that in this paper.

The case $t=0$ yields the desired criterion for the space $V/W$ to admit a symplectic resolution. It is a beautiful result due to Ginzburg--Kaledin \cite{gika} and Namikawa \cite{na} that $V/W$ admits a symplectic resolution if and only if the spectrum of the centre of ${\bf H}_{0,\bc}(W)$ is smooth for generic $\bc$. The space $X_{\bc}(W):={\rm Spec}(Z({\bf H}_{0,\bc}(W)))$ is called \emph{generalised Calogero--Moser space}. In \cite{go} Gordon introduced and studied extensively a finite-dimensional quotient of ${\bf H}_{0,\bc}(W)$,  called the  \emph{restricted rational Cherednik algebra}, whose simple modules are parametrised by $\Irr(W)$. The decomposition of this algebra into blocks induces a partition of $\Irr(W)$, known as \emph{Calogero--Moser partition}. We have that $X_{\bc}(W)$ is smooth if and only if the Calogero--Moser partition is trivial for all parabolic subgroups of $W$.  Following the classification of irreducible complex reflection groups, and the works of Etingof--Ginzburg \cite{EG}, Gordon \cite{go} and Gordon--Martino \cite{GM}, Bellamy \cite{bel} was able to prove that $V/W$ admits a symplectic resolution if and only if $W=G(\ell,1,n)$ or $W=G_4$.

It is believed that there exists a connection between the Calogero--Moser partition and the families of characters, first suggested by Gordon and Martino \cite{GM} for type $B_n$. In every case studied so far, the partition into Rouquier families (for a suitably chosen cyclotomic Hecke algebra) refines the Calogero--Moser partition (``Martino's conjecture''), while for finite Coxeter groups the two partitions coincide. The reasons for this connection are still unknown, since there is no apparent connection between Hecke algebras and rational Cherednik algebras at $t=0$. Inspired by this, and in an effort to construct a generalised Kazhdan--Lusztig cell theory, Bonnaf\'e and Rouquier have used the Calogero--Moser partition to develop a ``Calogero--Moser cell theory'' which can be applied to all complex reflection groups \cite{Bora}. The fruits of this very recent approach remain to be seen.

\subsection{Piece of notation and definition of blocks}\label{definitionofblocks}
Let $R$ be a commutative integral domain and let $F$ be  the field of fractions of $R$. Let $A$ be an $R$-algebra, free and finitely generated as an $R$-module. If $R'$ is a commutative integral domain containing $R$, we will write $R'A$ for $R' \otimes_R A$. 

Let now $K$ be a field containing $F$ such that the algebra $KA$ is semisimple. The primitive idempotents of the centre $Z(KA)$ of $KA$ are in bijection with the irreducible representations of $KA$. Let $\Irr(KA)$ denote the set of  irreducible representations of $KA$. For $\chi \in \Irr(KA)$, let $e_\chi$ be the corresponding primitive idempotent of $Z(KA)$. There exists a unique partition ${\rm Bl}(A)$ of $\Irr(KA)$ that is the finest with respect to the property:
$$ \forall B \in {\rm Bl}(A),\,\,\,\,\,\,e_B:=\sum_{\chi \in B} e_\chi \in A.$$
The elements $\{e_B\}_{B \in {\rm Bl}(A)}$ are the primitive idempotents of $Z(A)$. We have $A \cong \prod_{B \in {\rm Bl}(A)} Ae_B$. The parts of ${\rm Bl}(A)$ are the {\em blocks} of $A$.

\section{Iwahori--Hecke Algebras}
In this section we will focus on real reflection groups, while in the next section we will see what happens in the complex case.

\subsection{Kazhdan--Lusztig cells}\label{KLCells} Let $(W,S)$ be a finite Coxeter system. By definition, $W$ has a presentation of the form
$$W=\langle \,S\,|\, (st)^{m_{st}}=1\,\,\,\,\forall \, s,t \in S\, \rangle$$
with $m_{ss}=1$ and $m_{st} \geq 2$ for $s \neq t$. We have a  \emph{length function} $\ell: W \rightarrow \mathbb{Z}_{\geq 0}$  defined by $\ell(w):={\rm min}\,\{\,r\,|\,w=s_{i_1}\ldots s_{i_r}\,\,\text{with}\, \,s_{i_j} \in S\,\}$ for all $w \in W$.

Let $L:W \rightarrow\mathbb{Z}_{\geq 0}$ be a weight function, that is, a map such that $L(ww')=L(w)+L(w')$ whenever  $\ell(ww')=\ell(w)+\ell(w')$. For $s,\,t \in S$, we have $L(s)=L(t)$
whenever $s$ and $t$ are conjugate in $W$. 
Let $q$ be an indeterminate. We define the \emph{Iwahori--Hecke algebra}  of $W$ with parameter $L$ , denoted by $\mathcal{H}(W,L)$, to be the $\Z[q,q^{-1}]$-algebra generated by elements
$(T_s)_{s\in S}$ satisfying the relations:
$$(T_s-q^{L(s)})(T_s+q^{-L(s)})=0\,\,\,\,\,\,\,\text{and}\,\,\,\,\,\,\,
\underbrace{T_sT_tT_sT_t\ldots}_{m_{st}}=\underbrace{T_tT_sT_tT_s\ldots}_{m_{st}} \,\,\,\,\,\text{for }\, s\neq t.$$
If $L(s)=L(t)$ for all $s,t \in S$, we say that we are in the \emph{equal parameter case}. Since $L$ is a weight function, unequal parameters can only occur in irreducible types $B_n$, $F_4$ and dihedral groups $I_2(m)$ for $m$ even.

\begin{exmp} {\rm Let $W=\mathfrak{S}_3$. We have $W=\langle s, t\,|\, s^2=t^2=(st)^{3}=1 \rangle$. Let $l:=L(s)=L(t) \in \mathbb{Z}_{\geq 0}$. We have}
$$\mathcal{H}(W, l)= \langle T_s, T_t\,|\,T_sT_tT_s=T_tT_sT_t,\,  (T_s-q^l)(T_s+q^{-l})=(T_t-q^{l})(T_t+q^{-l})=0 \rangle.$$
\end{exmp}

Let $w\in W$ and let $w=s_{i_1}\ldots s_{i_r}$ be a reduced expression for $w$, that is, $r = \ell(w)$. Set $T_w:=T_{s_{i_1}}\ldots T_{s_{i_r}}$.
As a  $\Z[q,q^{-1}]$-module, $\mathcal{H}(W,L)$ is generated by the elements $(T_w)_{w\in W}$ satisfying the following multiplication formulas:
$$\left\{\begin{array}{ll}
T_{s}^2 =1 + (q^{L(s)}-q^{-L(s)})\,T_s& \text{ for }s\in S,\\& \\
T_wT_{w'} =T_{ww'}  & \text{ if } \ell(ww')=\ell(w)+\ell(w').
\end{array}\right.$$
The elements $(T_w)_{w\in W}$ form a basis of $\mathcal{H}(W,L)$, the \emph{standard basis}.

Let $i$ be the  algebra involution on $\mathcal{H}(W,L)$ given by $i(q)=q^{-1}$ and $i(T_s)=T_s^{-1}$ for $s \in S$ (as a consequence, $i(T_w)=T_{w^{-1}}^{-1}$ for all $w \in W$). 
By \cite[Theorem 1.1]{KL} (see \cite[Proposition 2]{219} for the unequal parameter case), for each $w \in W$, there exists an element $C_w \in \mathcal{H}(W,L)$ uniquely determined by the conditions
$$i(C_w)=C_w \,\,\,\,\,\text{ and }\,\,\,\,\,i(C_w)=T_w + \sum_{x \in W,\, x < w} P_{x,w}\,T_x\,,$$
where $<$ stands for the Chevalley--Bruhat order on $W$ and $P_{x,w} \in q^{-1}\mathbb{Z}[q^{-1}]$.
The elements $(C_w)_{w \in W}$   also form a basis of $\mathcal{H}(W,L)$, the {\em Kazhdan--Lusztig basis}.

\begin{exmp} {\rm  We have $C_1=T_1=1$ and,  for all $s \in S$, 
$C_s=\left\{\begin{array}{ll}
T_s & \text{ if } L(s)=0 \\
T_s+q^{-L(s)}T_1 & \text{ if } L(s)>0
\end{array}\right.$.}
\end{exmp}

Using the Kazhdan--Lusztig basis, we can now define the following three preorders on $W$. For $x,\,y \in W$, we have
\begin{itemize}
\item $x \leq_{\mathcal{L}} y$\, if $C_x$ appears with non-zero coefficient  in $hC_y$ for some $h \in \mathcal{H}(W,L)$. \smallbreak
\item $x\leq_{\mathcal{R}} y$\, if $C_x$ appears with non-zero coefficient  in $C_yh'$ for some $h' \in \mathcal{H}(W,L)$. \smallbreak
\item $x \leq_{\mathcal{LR}} y$\, if $C_x$ appears with non-zero coefficient  in $hC_yh'$ for some $h,\,h' \in \mathcal{H}(W,L)$.
\end{itemize}

The preorder $\leq_{\mathcal{L}} $ defines an equivalence relation $\sim_{\mathcal{L}}$ on $W$ as follows:
$$x \sim_{\mathcal{L}}  y \,\,\Leftrightarrow\,\, x \leq_{\mathcal{L}}  y \,\text{ and }\,  y \leq_{\mathcal{L}}  x.$$
The equivalence classes for $\sim_{\mathcal{L}}$ are called {\em left cells}. 
Similarly, one can define equivalence relations $\sim_{\mathcal{R}}$ and  $\sim_{\mathcal{LR}}$ on $W$, whose equivalence classes are called, respectively, {\em right cells} and {\em two-sided cells}.

\begin{exmp} {\rm For $W=\mathfrak{S}_3=\{1,s,t,st,ts,sts=tst\}$ and $l>0$, 
\begin{itemize}
\item the left cells are $\{1\}$,\,$\{s,ts\}$,\,$\{t,st\}$ and $\{sts\}$ ; \smallbreak
\item the right cells are $\{1\}$,\,$\{s,st\}$,\,$\{t,ts\}$ and $\{sts\}$ ; \smallbreak
\item the two-sided cells are $\{1\}$,\,$\{s,t,st,ts\}$ and $\{sts\}$.
\end{itemize}
If $l=0$, then all elements of $W$ belong to the same cell (left, right or two-sided).}
\end{exmp}

Let now $\mathfrak{C}$ be a left cell of $W$. The following two $\Z[q,q^{-1}]$-modules are left ideals  of $\mathcal{H}(W,L)$:
$$\mathcal{H}_{\leq_{\mathcal{L}}\mathfrak{C}} = \langle C_y\,|\,y \leq_{\mathcal{L}}w, w \in \mathfrak{C} \rangle_{\Z[q,q^{-1}]}\,\,\,\,\,\,\,\,\text{ and }\,\,\,\,\,\,\,\,
\mathcal{H}_{<_{\mathcal{L}}\mathfrak{C}} = \langle C_y\,|\,y \leq_{\mathcal{L}}w, w \in \mathfrak{C} , y \notin \mathfrak{C} \rangle_{\Z[q,q^{-1}]}.$$
Then
$$\mathcal{M}_\mathfrak{C}:= \mathcal{H}_{\leq_{\mathcal{L}}\mathfrak{C}} / \mathcal{H}_{<_{\mathcal{L}}\mathfrak{C}}$$ is a free  left $\mathcal{H}(W,L)$-module
with basis  indexed by the elements of $\mathfrak{C}$.

Let $K$ be a field containing $\mathbb{Z}[q,q^{-1}]$ such that the algebra $K\mathcal{H}(W,L)$ is split semisimple (for example, take $K=\mathbb{C}(q)$). Then, since the left cells form a partition of $W$, we obtain a corresponding direct sum decomposition of $K\mathcal{H}(W,L)$:
\begin{equation}\label{celldecomp}
K\mathcal{H}(W,L) \cong \bigoplus_{\mathfrak{C} \text{ left cell }} K\mathcal{M}_\mathfrak{C} \quad \quad \text{(isomorphism of left $K\mathcal{H}(W,L)$-modules),}
\end{equation}
where $K\mathcal{M}_\mathfrak{C}:=K \otimes_{\mathbb{Z}[q,q^{-1}]}\mathcal{M}_\mathfrak{C}$.
We obtain analogous decompositions with respect to right and two-sided cells.

\subsection{Schur elements and the $a$-function}\label{schurreal}
 From now on, set $R:= \Z[q,q^{-1}]$ and let $K$ be a field containing $R$ such that the algebra $K\mathcal{H}(W,L)$ is split semisimple. 

Using the standard basis of the Iwahori--Hecke algebra, we define the linear map $\tau:\mathcal{H}(W,L)\rightarrow R$  by setting 
$$\tau(T_w):=\left\{\begin{array}{ll}
1& \text{ if }w=1,\\
0& \text{ otherwise}. 
\end{array}\right.$$
The map $\tau$ is a symmetrising trace on $\mathcal{H}(W,L)$, that is, 
\begin{enumerate}[(a)]
\item $\tau(hh')=\tau(h'h)$ for all $h,h' \in \mathcal{H}(W,L)$, and \smallbreak
\item the map $\widehat{\tau}: \mathcal{H}(W,L) \rightarrow \mathrm{Hom}_R(\mathcal{H}(W,L),R)$, $h \mapsto (x \mapsto \tau(hx))$ is an isomorphism of $\mathcal{H}(W,L)$-bimodules.
\end{enumerate}
Moreover, the elements $(T_{w^{-1}})_{w \in W}$ form a basis of $\mathcal{H}(W,L)$ dual to the standard basis with respect to $\tau$ (that is, $\tau(T_{w^{-1}} T_{w'})=\delta_{w,w'}$) \cite[Proposition 8.1.1]{gepf}.
The map $\tau$ is called the {\em canonical symmetrising trace} on $\mathcal{H}(W,L)$, because it specialises to the canonical symmetrising trace on the group algebra $\Z[W]$ when $q \mapsto 1$. 

Now, the map $\tau$ can be extended to $K\mathcal{H}(W,L)$ by extension of scalars. By Tits's deformation theorem (see, for example, \cite[Theorem 7.4.6]{gepf}), the specialisation $q \mapsto 1$ induces a bijection between the set of irreducible representations $\Irr(K\mathcal{H}(W,L))$ of $K\mathcal{H}(W,L)$ and   the set of irreducible representations $\Irr(W)$ of $W$. 
For $E \in \Irr(W)$, let $\chi_E$ be the corresponding irreducible character of  $K\mathcal{H}(W,L)$
and let $\omega_{\chi_E}$ be the corresponding central character. We define
$$s_E:= \chi_E(\widehat{\tau}^{-1}(\chi_E))/\chi_E(1) = \omega_{\chi_E}(\widehat{\tau}^{-1}(\chi_E))$$
to be the \emph{Schur element} of $\mathcal{H}(W,L)$ associated with $E$.
Geck has shown (see \cite[Proposition 7.3.9]{gepf}) that $s_E \in \Z_K[q,q^{-1}]$ for all $E \in \Irr(W)$, where $\Z_K$ denotes the integral closure of $\Z$ in $K$. 
We have
\begin{equation}
\tau = \sum_{E \in \Irr(W)}\frac{1}{s_E}\chi_E
\end{equation}
and 
\begin{equation}\label{symmetricblock}
e_E = \frac{1}{s_E} \sum_{w \in W} \chi_E(T_w)\, T_{w^{-1}}\,,
\end{equation}
where $e_E$ is the primitive central idempotent of $K\mathcal{H}(W,L)$  corresponding to $E$.
Both results are due to Curtis and Reiner \cite{CuRe}, but we follow the exposition in  
 \cite[Theorem 7.2.6]{gepf} and \cite[Proposition 7.2.7]{gepf} respectively.

\begin{exmp} {\rm In the group algebra case ($L(s)=0$ for all $s\in S$), we have $s_E=|W|/\chi_E(1)$ for all $E \in \Irr(W)$.  }
\end{exmp}

\begin{exmp} {\rm The irreducible representations of the symmetric group $\mathfrak{S}_n$ are parametrised by the partitions of $n$.
For $W=\mathfrak{S}_3$, there are three irreducible representations.
Let $E^{(3)}$, $E^{(2,1)}$ and $E^{(1,1,1)}$ denote respectively the trivial, reflection and sign representation of $\mathfrak{S}_3$. We have
\begin{center}
$s_{E^{(3)}}=(q^{2l}+1)(q^{4l}+q^{2l}+1),\,s_{E^{(2,1)}}=q^{2l}+1+q^{-2l},\,s_{E^{(1,1,1)}}=(q^{-2l}+1)(q^{-4l}+q^{-2l}+1).$
\end{center}}
\end{exmp}

We can define the functions $a: \Irr(W) \rightarrow \Z$ and $A: \Irr(W) \rightarrow \Z$ by setting
$$a(E):=-{\rm val}_q(s_E) \,\,\,\,\,\text{ and }\,\,\,\,\,A(E):=-{\rm deg}_q(s_E).$$
Note that both functions depend on $L$. For brevity, we will write $a_E$ for $a(E)$ and $A_E$ for $A(E)$.

\begin{exmp} {\rm For $W=\mathfrak{S}_3$, we have
\begin{center}
$a_{E^{(3)}}=0,\, a_{E^{(2,1)}}=2l,\,a_{E^{(1,1,1)}}=6l$ \,\,\,\,\,and\,\,\,\,\,
$A_{E^{(3)}}=-6l,\,A_{E^{(2,1)}}=-2l,\,A_{E^{(1,1,1)}}=0.$
\end{center}}
\end{exmp}

The Schur elements of $\mathcal{H}(W,L)$ have been explicitly calculated for all finite Coxeter groups: 
\begin{itemize}
\item for type $A_n$ by Steinberg \cite{St1}, \smallbreak
\item for type $B_n$ by Hoefsmit \cite{Ho}, \smallbreak
\item for type  $D_n$ by Benson and Gay \cite{BeGa} (it derives from type $B_n$ with the use of Clifford theory), \smallbreak
\item for dihedral groups $I_2(m)$ by Kilmoyer and Solomon \cite{Kiso}, \smallbreak
\item for $F_4$ by Lusztig \cite{Lu79b}, \smallbreak
  \item for $E_6$ and $E_7$ by Surowski \cite{Sur78}, \smallbreak
  \item for $E_8$ by Benson \cite{Ben79},\smallbreak
  \item for $H_3$ by Lusztig  \cite{Lu82},\smallbreak
  \item for $H_4$ by Alvis and Lusztig \cite{AlLu82}.
\end{itemize}
There have been other subsequent proofs of the above results. For example, Iwahori--Hecke algebras of types $A_n$ and $B_n$ are special cases of Ariki--Koike algebras, whose Schur elements have been independently obtained by Geck--Iancu--Malle \cite{GIM00} and Mathas \cite{Mathas}.

A case-by-case analysis shows that the Schur elements of $\mathcal{H}(W,L)$ can be written in the form
\begin{equation}\label{cycox}
s_E = \xi_E \,q^{-a_E} \prod_{\Phi \in {\rm Cyc}_E} \Phi(q^{n_{E,\Phi}})
\end{equation}
where $\xi_E \in \Z_K$, $n_{E,\Phi} \in \Z_{>0}$ and ${\rm Cyc}_E$ is a family of $K$-cyclotomic polynomials (see \cite[Chapters 10 and 11]{gepf}, \cite[Theorem 4.2.5]{mybook}).

\begin{exmp} {\rm For $W=\mathfrak{S}_3$, we have
\begin{center}
$s_{E^{(3)}}=\Phi_2(q^{2l})\Phi_3(q^{2l}),\,s_{E^{(2,1)}}=q^{-2l}\Phi_3(q^{2l}),\,s_{E^{(1,1,1)}}=q^{-6l}\Phi_2(q^{2l})\Phi_3(q^{2l}).$
\end{center}}
\end{exmp}

\subsection{Families of characters and Rouquier families}\label{fam} 

The {\em families of characters} are a special partition of the set of irreducible representations of $W$.  In the case where $W$ is a Weyl group, these families play an essential role in the definition of the families of unipotent characters for the corresponding finite reductive groups. Their original definition is due to Lusztig \cite[4.2]{26} and uses the $a$-function.

Let $I \subseteq S$ and consider the parabolic subgroup $W_I \subseteq W$ generated by $I$. Then we have a corresponding parabolic subalgebra $\mathcal{H}(W_I,L) \subseteq \mathcal{H}(W,L)$. By extension of scalars from $R$ to $K$, we also have a subalgebra $K\mathcal{H}(W_I,L) \subseteq K\mathcal{H}(W,L)$,  and a corresponding $a$-function on the set of irreducible representations of $W_I$. Denote by ${\rm Ind}^S_I$ the induction of representations from $W_I$ to $W$. Let $E \in \Irr(W)$ and $M \in \Irr(W_I)$. We will write 
$M \leadsto_L E$ if $E$ is a constituent of ${\rm Ind}^S_I (M)$ and $a_E = a_M$.

\begin{defn}{\rm 
The partition of $\Irr(W)$ into \emph{families} is
defined inductively as follows: When $W = \{1\}$, there is only one family; it consists of the
unit representation of $W$. Assume now that $W \neq \{1\}$ and that the families have
already been defined for all proper parabolic subgroups of $W$. Then $E,E' \in
\Irr(W)$ are  in the same family of $W$ if there exists a finite sequence
$E = E_0,E_1,\ldots,E_r = E'$ in $\Irr(W)$ such that, for each $i \in  \{0,1,...,r - 1\}$,
the following condition is satisfied: There exist a subset $I_i \subsetneqq S$ and $M_i,\,M_i' \in
\Irr(W_{I_i})$ such that $M_i,M_i'$ belong to the same family of $W_{I_i}$  and either 
$$M_i \leadsto_L  E_i \,\,\,\,\text{ and }\,\,\,\, M_i' \leadsto_L  E_{i+1}$$
or
$$M_i \leadsto_L  E_i \otimes \varepsilon \,\,\,\,\text{ and }\,\,\,\, M_i' \leadsto_L  E_{i+1} \otimes \varepsilon,$$
where $\varepsilon$ denotes the sign representation of $W$. We will also refer to these families as \emph{Lusztig families}.}
\end{defn}

Lusztig \cite[3.3 and 3.4]{lu3} has shown that the functions $a$ and $A$ are both constant on the families of characters, that is, if $E$ and $E'$ belong to the same family, then $a_E=a_{E'}$ and $A_E=A_{E'}$.

The decomposition of $W$ into two-sided cells can be used to facilitate the description of the partition of $\Irr(W)$ into families of characters.
As we saw in the previous subsection, Tits's deformation theorem yields a bijection between $\Irr(K\mathcal{H}(W,L))$ and  $\Irr(W)$. Let $E \in \Irr(W)$ and let $V^E$ be the corresponding simple module of $K\mathcal{H}(W,L)$. Following the direct sum decomposition given by $(\ref{celldecomp})$, there exists a left cell $\mathfrak{C}$ such that $V^E$ is a constituent of $\mathcal{M}_\mathfrak{C}$; furthermore, all such left cells are contained in the same two-sided cell. This two-sided cell, therefore, only depends on $E$ and will be denoted by $\mathcal{F}_E$. Thus, we obtain a natural surjective map
$$\Irr(W) \rightarrow \{\text{set of two-sided cells of } W\},\,\, E \mapsto \mathcal{F}_E$$
(see \cite[5.15] {26} for the equal parameter case; the same argument works in general). 

\begin{defn} 
{\rm Let $E,E' \in \Irr(W)$. We will say that $E$ and $E'$ belong to the same \emph{Kazhdan--Lusztig family} if $\mathcal{F}_E=\mathcal{F}_{E'}$.}
\end{defn}

The following remarkable result, relating Lusztig families and Kazhdan--Lusztig families, has been proved by Barbasch--Vogan and Lusztig  for finite Weyl groups in the equal parameter case \cite[5.25]{26}. It was subsequently proved, in \cite[23.3]{31} and \cite{13}, to hold for any finite Coxeter group and  any weight function $L$, assuming that Lusztig's conjectures {\bf P1--P15} \cite[14.2]{31} are satisfied.

\begin{thm}
Assume that Lusztig's conjectures {\bf P1--P15} hold. The Lusztig families and the Kazhdan--Lusztig families coincide. 
\end{thm}

Lusztig's conjectures  {\bf P1--P15} concern properties of the Kazhdan--Lusztig basis which should hold for any Coxeter group and in the general multiparameter case. 
For the moment, Conjectures {\bf P1--P15} have been proved in the following cases: 
\begin{itemize}
\item Equal parameter case for finite Weyl groups  \cite{KLautre,31,springer}. \smallbreak
\item Equal parameter case for $H_3$, $H_4$ and dihedral groups $I_2(m)$  \cite{Al, DuCloux}. \smallbreak
\item Unequal parameter case for $F_4$ and dihedral groups $I_2(m)$ \cite{Gecells, GeP1}. \smallbreak
\item Asymptotic case and some other cases for $B_n$ \cite{BI, Bo, BGIL}.
\end{itemize}
Moreover, these are exactly the cases where we have a description of the Kazhdan--Lusztig cells and Kazhdan--Lusztig families. A conjectural combinatorial 
description of the Kazhdan--Lusztig cells for type $B_n$ is given by \cite{BGIL}. 

\begin{exmp} {\rm The group $\mathfrak{S}_3$ has three irreducible representations.
For $l>0$, each irreducible representation forms a family on its own. This is true in general for the symmetric group $\mathfrak{S}_n$. For $l=0$, all irreducible representations belong to the same family. This is true in general for the group algebra ($L(s)=0$ for all $s\in S$) of every finite Coxeter group.}
\end{exmp}

In \cite{Ro} Rouquier gave an alternative definition for Lusztig's families. He showed that, for finite Weyl groups in the equal parameter case, the families of characters coincide with the blocks of the Iwahori--Hecke algebra $\mathcal{H}(W,L)$  over the {\em Rouquier ring}
$$\mathcal{R}_K(q):=\Z_K[q,q^{-1},(q^n-1)^{-1}_{n \geq1}] ,$$
that is, following (\ref{symmetricblock}), the non-empty subsets $B$ of $\Irr(W)$ which are minimal with respect to the property:
$$ \sum_{E \in B} \frac{\chi_E(h)}{s_E} \in \mathcal{R}_K(q)\,\,\,\,\,\forall \,h \in \mathcal{H}(W,L).$$
These are the {\em Rouquier families} of $\mathcal{H}(W,L)$. One advantage of this definition, as we will see in the next section, is that it can be also applied to complex reflection groups. This is important in the project ``Spetses'' \cite{BMM2, BMM3}.

Following the determination of Rouquier families for all complex reflection groups (see \S\ref{Roublo} for references), and thus for all finite Coxeter groups, one can check that Rouquier's result holds for all finite Coxeter groups for all choices of parameters (by comparing the Rouquier families with the already known  Lusztig families \cite{26,31}); that is, we have the following:

\begin{thm}
Let $(W,S)$ be a finite Coxeter system and let $\mathcal{H}(W,L)$ be an Iwahori--Hecke algebra associated to $W$. The Lusztig families and the Rouquier families of $\mathcal{H}(W,L)$ coincide. 
\end{thm}

\subsection{Canonical basic sets}\label{canbasset} 
As we saw in \S\ref{fam}, 
the specialisation $q \mapsto 1$ yields a bijection between the set of irreducible representations of $K\mathcal{H}(W,L)$ and $\Irr(W)$. What happens though when $q$ specialises to a complex number? The resulting Iwahori--Hecke algebra is not necessarily semisimple and the first questions that need to be answered are the following: What are the simple modules for the newly obtained algebra? Is there a good way to parametrise them? What are their dimensions? One major approach to answering these questions is through the existence of ``canonical basic sets''. 

Let $\theta: \Z_K[q,q^{-1}] \rightarrow K(\eta),\, q \mapsto \eta$ be a ring homomorphism such that $\eta$ is a non-zero complex number. Let us denote by $\mathcal{H}_\eta(W,L)$ the algebra obtained as a specialisation of $\mathcal{H}(W,L)$ via $\theta$. Set $\K:=K(\eta)$. We have the following semisimplicity criterion \cite[Theorem 7.4.7]{gepf}:

\begin{thm}\label{sscriterion} The algebra $\K\mathcal{H}_\eta(W,L)$ is semisimple if and only $\theta(s_E) \neq 0$ for all $E \in \Irr(W)$.
\end{thm}

Following  (\ref{cycox}), $\K\mathcal{H}_\eta(W,L)$ is semisimple unless $\eta$ is a root of unity. 

\begin{exmp}
{\rm The algebra $\Q(\eta)\mathcal{H}_\eta(\mathfrak{S}_3,l)$ is semisimple if and only if $\eta^{2l} \notin \{-1,\omega,\omega^2\}$, where $\omega:={\rm exp}(2\pi i/3)$.}
\end{exmp}

If $\K\mathcal{H}_\eta(W,L)$ is semisimple, then, by Tits's deformation theorem, the specialisation $\theta$ yields a bijection between $\Irr(K\mathcal{H}(W,L))$ and $\Irr(\K\mathcal{H}_\eta(W,L))$. Thus, the irreducible representations of $\K\mathcal{H}_\eta(W,L)$ are parametrised by $\Irr(W)$. Hence, we need to see what happens when $\K\mathcal{H}_\eta(W,L)$ is not semisimple.

  Let 
$R_0(K\mathcal{H}(W,L))$ (respectively $R_0(\K\mathcal{H}_\eta(W,L))$) be the Grothendieck group
 of finitely generated $K\mathcal{H}(W,L)$-modules (respectively $\K\mathcal{H}_\eta(W,L)$-modules). 
  It is generated by the classes $[U]$ of the simple  $K\mathcal{H}(W,L)$-modules (respectively $\K\mathcal{H}_\eta(W,L)$-modules) $U$.
   Then we obtain a well-defined  decomposition map
$$d_{\theta}:R_0 (K\mathcal{H}(W,L))    \to R_0 (\K\mathcal{H}_\eta(W,L))$$
such that, for all $E \in \Irr(W)$, we have
$$d_{\theta}([V^E])=\sum_{M\in \Irr(\K\mathcal{H}_\eta(W,L))}
 [ V^E : M] [M].$$
The matrix
\[ D_{\theta}=\left([ V^E :M]\right)_{E \in \Irr(W),\,M \in \Irr(\K\mathcal{H}_\eta(W,L)) }\]
is called the \emph{decomposition matrix with respect to} $\theta$. If $\K\mathcal{H}_\eta(W,L)$ is semisimple, then $D_\theta$ is  a permutation matrix.

\begin{defn}{\rm A {\em canonical basic set} with respect to $\theta$ is a subset $\mathcal{B}_\theta$ of $\Irr(W)$ such that
\begin{enumerate}[(a)]
\item there exists a bijection  $\Irr(\K\mathcal{H}_\eta(W,L)) \rightarrow \mathcal{B}_\theta$, $M \mapsto E_M$; \smallbreak
\item $[ V^{E_M} :M] =1$ for all $M \in \Irr(\K\mathcal{H}_\eta(W,L))$; \smallbreak
\item if $[ V^E :M] \neq 0$ for some $E \in \Irr(W),\,M \in \Irr(\K\mathcal{H}_\eta(W,L))$, then either $E=E_M$ or $ a_{E_M}<a_E $.
\end{enumerate}
}
\end{defn}

If a canonical basic set exists, the decomposition matrix has a lower unitriangular form (with an appropriate ordering of the rows).
Thus, we can obtain a lot of information about the simple modules of $\K\mathcal{H}_\eta(W,L)$ 
from what we already know about the simple modules of  $K\mathcal{H}(W,L)$.

A general existence result for canonical basic sets is proved by Geck in \cite[Theorem 6.6]{Gsurvey}, following earlier work of Geck \cite{96}, Geck--Rouquier \cite{GR} and Geck--Jacon \cite{GJ}. Another proof is given in \cite{GJlivre}. In every case canonical basic sets are known explicitly, thanks to the work of many people. For a complete survey on the topic, we refer the reader to \cite{GJlivre}.

\begin{exmp}\label{excanbassym}{\rm Let $W$ be the symmetric group $\mathfrak{S}_n$. Then $W$ is generated by the transpositions $s_i=(i,i+1)$ for all $i=1,\ldots,n-1$, which are all conjugate in $W$. Set $l:=L(s_1)$  
and let $\eta^{2l}$ be a primitive root of unity of order $e>1$. By \cite[Theorem 7.6]{DJ0}, we have that, in this case, the canonical basic set $\mathcal{B}_\theta$  is the set of $e$-regular partitions (a partition is $e$-regular if it does not have $e$ non-zero equal parts). For example, for $n=3$, we have $\mathcal{B}_\theta=\{E^{(3)}, E^{(2,1)}\}$ for $e \in \{2,3\}$, and  $\mathcal{B}_\theta=\Irr(\mathfrak{S}_3)$ for $e>3$.
}
\end{exmp}

\section{Cyclotomic Hecke Algebras}

Cyclotomic Hecke algebras generalise the notion of Iwahori--Hecke algebras to the case of complex reflection groups.
For any positive integer $e$ we will write $\zeta_e$ for $\exp(2\pi i/e)\in \C$.

\subsection{Hecke algebras for complex reflection groups}\label{Hecke}

Let $\mathfrak{h}$ be a finite dimensional complex vector space. A \emph{pseudo-reflection} is a non-trivial element $s \in  \mathrm{GL}(\mathfrak{h})$ that fixes a hyperplane pointwise,
that is, ${\rm dim}_{\C}{\rm Ker}(s - {\rm id}_\mathfrak{h})={\rm dim}_{\C}\mathfrak{h}-1$. The hyperplane ${\rm Ker}(s - {\rm id}_\mathfrak{h})$ is the \emph{reflecting hyperplane} of $s$.
A {\em complex reflection group} is a finite subgroup of $\mathrm{GL}(\mathfrak{h})$ generated by pseudo-reflections. The classification of (irreducible) complex reflection groups is due to Shephard and Todd \cite{ShTo}:

\begin{thm}\label{ShToClas} Let $W \subset \mathrm{GL}(\mathfrak{h})$ be an irreducible complex
reflection group (i.e., $W$ acts irreducibly on $\mathfrak{h}$). Then one of
the following assertions is true:
\begin{itemize}
  \item There exist positive integers $\ell,p,n$ with $\ell/p \in \Z$ and $\ell > 1$ such
  that $(W,\mathfrak{h}) \cong (G(\ell,p,n),\C^n)$, where $G(\ell,p,n)$ is the group of all 
  $n \times n$ monomial matrices whose non-zero entries are ${\ell}$-th roots of unity, while the product of all non-zero
  entries is an $(\ell/p)$-th root of unity. \smallbreak
  \item There exists a positive integer $n$ such that $(W,\mathfrak{h}) \cong (\mathfrak{S}_n, \C^{n-1})$.\smallbreak
  \item $(W,\mathfrak{h})$ is isomorphic to one of the 34 exceptional groups
  $G_n$ $(n=4,\ldots,37)$.
\end{itemize}
\end{thm}

\begin{rem}
{\rm We have $G(1,1,n) \cong \mathfrak{S}_n$, $G(2,1,n) \cong B_n$,  $G(2,2,n) \cong D_n$, $G(m,m,2) \cong I_2(m)$, \\
$G_{23} \cong H_3$,  $G_{28}  \cong  F_4$, $G_{30}  \cong H_4$, $G_{35}  \cong  E_6$, $G_{36}  \cong  E_7$, $G_{37}  \cong E_8$.       
}
\end{rem}

Let $W \subset \mathrm{GL}(\mathfrak{h})$ be a complex
reflection group. Benard \cite{Ben} and Bessis \cite{Bes1} have proved (using a case-by-case
analysis) that the field $K$ generated by the traces on $\mathfrak{h}$ of all the elements of $W$ is a splitting field for $W$. The field $K$  is called the \emph{field of
definition} of $W$.  If $K \subseteq \mathbb{R}$, then $W$ is a finite Coxeter group, and
 if $K=\mathbb{Q}$, then $W$ is a Weyl group.

Let $\mathcal{A}$ be the set of reflecting hyperplanes of $W$.
Let $\mathfrak{h}^{\textrm{reg}}:= \mathfrak{h}\setminus \bigcup_{H\in \mathcal{A}} H$ and $B_W := \pi_1(\mathfrak{h}^{\textrm{reg}}/W, x_0)$, where $x_0$ is some fixed basepoint.
The group $B_W$ is the {\em braid group} of $W$.
For every orbit $\mathcal{C}$ of $W$ on $\mathcal{A}$, we set
$e_{\mathcal{C}}$ the common order of the subgroups $W_H$, where $H$
is any element of $\mathcal{C}$ and $W_H$ is the pointwise stabiliser of $H$. Note that $W_H$ is cyclic, for all $H \in \mathcal{A}$.

We choose a set of indeterminates
$\textbf{u}=(u_{\mathcal{C},j})_{(\mathcal{C} \in
\mathcal{A}/W)(0\leq j \leq e_{\mathcal{C}}-1)}$ and we denote by
$\mathbb{Z}[\textbf{u},\textbf{u}^{-1}]$ the Laurent polynomial ring
in all the indeterminates $\textbf{u}$. We define the \emph{generic
Hecke algebra} \index{generic Hecke algebra} $\mathcal{H}(W)$ of $W$ to be the quotient of the group
algebra $\mathbb{Z}[\textbf{u},\textbf{u}^{-1}]B_W$ by the ideal
generated by the elements of the form
$$(\textbf{s}-u_{\mathcal{C},0})(\textbf{s}-u_{\mathcal{C},1}) \cdots (\textbf{s}-u_{\mathcal{C},e_{\mathcal{C}}-1}),$$
where $\mathcal{C}$ runs over the set $\mathcal{A}/W$ and
$\textbf{s}$ runs over the set of monodromy generators around 
the images in $\mathfrak{h}^{\textrm{reg}}/W$ of 
the elements of $\mathcal{C}$ \cite[\S 4]{BMR}.

From now on, we  will make certain assumptions for $\mathcal{H}(W)$.  These assumptions are known to hold 
for all finite Coxeter groups \cite[IV, \S 2]{Bou05}, $G(\ell,p,n)$ \cite{BMM2, MM98, GIM00} and a few of the exceptional complex reflection groups \cite{Marin1, Marin2}\footnote{In \cite{MM10} it is mentioned that these assumptions have been confirmed computationally by M\"uller  in several exceptional cases, but 
this work is not published. Moreover, in \cite{BM} the assumption (a) is proved for the groups $G_4$, $G_5$, $G_{12}$ and $G_{25}$, but Marin pointed out in \cite{Marin2} that these proofs might contain a questionable argument.}; they are expected to be true for all complex reflection groups.

\begin{hyp}\label{assumptions}\begin{enumerate}[(a)]
\item The algebra $\mathcal{H}(W)$ is a free {\rm $\mathbb{Z}[\textbf{u},\textbf{u}^{-1}]$}-module of rank equal to the order of $W$.\smallbreak
\item There exists a symmetrising trace $\tau$ on $\mathcal{H}(W)$ that satisfies certain canonicality conditions  \cite[\S1 and 2]{BMM2}; the form $\tau$ specialises to the canonical symmetrising form on the group algebra when  $u_{\mathcal{C},j}\mapsto \zeta_{e_\mathcal{C}}^{j}$.
\end{enumerate}
\end{hyp}

Under these assumptions, Malle \cite[5.2]{Ma4} has shown that there exists $N_W \in \Z_{>0}$ such that  if we take
\begin{equation}\label{NW}
u_{\mathcal{C},j} = \zeta_{e_\mathcal{C}}^{j}v_{\mathcal{C},j}^{N_W}
\end{equation}
and set $\textbf{{v}}:=(v_{\mathcal{C},j})_{(\mathcal{C} \in
\mathcal{A}/W)(0\leq j \leq e_{\mathcal{C}}-1)}$, then the 
$K(\textbf{{v}})$-algebra
$K(\textbf{{v}})\mathcal{H}(W)$ is split semisimple.
By Tits's deformation theorem, it follows
that the specialisation $v_{\mathcal{C},j}\mapsto 1$ induces a
bijection between
$\mathrm{Irr}(K(\textbf{v})\mathcal{H}(W))$ and  $\mathrm{Irr}(W)$. From now on, we will consider $\mathcal{H}(W)$ as an algebra over $\mathbb{Z}_K[\textbf{v},\textbf{v}^{-1}]$, where $\Z_K$ denotes the integral closure of $\Z$ in $K$.

\begin{exmp} {\rm The group $W=G(\ell,1,n)$ is isomorphic to the wreath product $(\Z/\ell \Z)\wr \mathfrak{S}_n$ and its splitting field is $K=\Q(\zeta_{\ell})$.  In this particular case, we can take $N_{W}=1$. The algebra $K(\textbf{v})\mathcal{H}(W)$ is generated by elements $\textbf{s}, \textbf{t}_1,\ldots,\textbf{t}_{n-1}$ satisfying the braid relations of type $B_n$, 
$$\textbf{s}\textbf{t}_1\textbf{s}\textbf{t}_1= \textbf{t}_1\textbf{s}\textbf{t}_1\textbf{s}, \,\,\,\,\textbf{s}\textbf{t}_i= \textbf{t}_i\textbf{s} \,\,\,\,\text{and}\,\,\,\,
\textbf{t}_{i-1}\textbf{t}_i\textbf{t}_{i-1}=\textbf{t}_i\textbf{t}_{i-1}\textbf{t}_i  \,\,\,\,\text {for } i=2,\ldots,n-1 ,\,\,\,\, \textbf{t}_{i}\textbf{t}_j= \textbf{t}_{j}\textbf{t}_i \,\,\,\,\text {for } |i-j|>1,$$
together with the extra relations
$$(\textbf{s}-v_{\textbf{s},0})(\textbf{s}-\zeta_\ell v_{\textbf{s},1})\cdots(\textbf{s}-\zeta_\ell^{\ell-1} v_{\textbf{s},\ell-1})=0\,\,\,\,\text{and}\,\,\,\,
(\textbf{t}_i-v_{\textbf{t},0})(\textbf{t}_i+ v_{\textbf{t},1})=0 \,\,\,\,\text {for all } i=1,\ldots,n-1.$$
The Hecke algebra of $G(\ell,1,n)$ is also known as \emph{Ariki--Koike algebra}, with the last quadratic relation usually looking like this:
$$({\textbf{t}}_i-q)({\textbf{t}}_i+ 1)=0,$$
where $q$ is an indeterminate.
The irreducible representations of $G(\ell,1,n)$, and thus the irreducible representations of  $K(\textbf{v})\mathcal{H}(W)$, are parametrised by the $\ell$-partitions of $n$.
}
\end{exmp}

Let now $q$ be an indeterminate and let $\textbf{m}=(m_{\mathcal{C},j})_{(\mathcal{C} \in
\mathcal{A}/W)(0\leq j \leq e_{\mathcal{C}}-1)}$ be a family of integers. The $\Z_K$-algebra morphism $${\varphi_{\textbf{m}}}: \mathbb{Z}_K[\textbf{{v}},\textbf{{v}}^{-1}] \rightarrow \mathbb{Z}_K[q,q^{-1}],\,\,v_{\mathcal{C},j} \mapsto q^{m_{\mathcal{C},j}}$$ is called a \emph{cyclotomic specialisation}. The $\mathbb{Z}_K[q,q^{-1}]$-algebra $\mathcal{H}_{\varphi_{\textbf{m}}}(W)$ obtained as the specialisation of $\mathcal{H}(W)$ via ${\varphi_{\textbf{m}}}$ is called a \emph{cyclotomic Hecke algebra} associated with $W$. The Iwahori--Hecke algebras defined in the previous section are cyclotomic Hecke algebras associated with real reflection groups. The algebra
$K(q)\mathcal{H}_{\varphi_{\textbf{m}}}(W)$ is split semisimple \cite[Proposition 4.3.4]{mybook}.
By Tits's deformation theorem,  the specialisation $q \mapsto 1$ yields a
bijection between
$\mathrm{Irr}(K(q)\mathcal{H}_{\varphi_{\textbf{m}}}(W))$ and  $\mathrm{Irr}(W)$.

\subsection{Schur elements and the $a$-function}

The symmetrising trace $\tau$ (see Hypothesis \ref{assumptions}) can be extended to $K(\textbf{v})\mathcal{H}(W)$ by extension of scalars, and can be used to define Schur elements
$(s_E)_{E \in \Irr(W)}$ for 
$\mathcal{H}(W)$. The Schur elements of $\mathcal{H}(W)$ have been explicitly calculated for all complex reflection groups:
\begin{itemize}
\item for finite Coxeter groups see \S \ref{schurreal} ;\smallbreak
\item for complex reflection groups of type $G(\ell,1,n)$ by Geck--Iancu--Malle \cite{GIM00} and Mathas \cite{Mathas}; \smallbreak
\item for complex reflection groups of type $G(\ell,2,2)$ by Malle \cite{Ma2}; \smallbreak
\item for the remaining exceptional complex reflection groups by Malle \cite{Ma2, Ma5}.\smallbreak
\end{itemize} 
With the use of Clifford theory, we obtain the Schur elements for type $G(\ell,p,n)$ from those of type $G(\ell,1,n)$ when $n>2$ or $n=2$ and $p$ is odd.
The Schur elements for type $G(\ell,p,2)$ when $p$ is even derive from those of type $G(\ell,2,2)$. See \cite{MaUni}, \cite[A.7]{mybook}.

Using a case-by-case analysis, we have been able to determine that the Schur elements of $\mathcal{H}(W)$  have the following form \cite[Theorem 4.2.5]{mybook}.

\begin{thm}\label{Schur element generic}
Let $E \in \Irr(W)$. The Schur element $s_E$  is an element of
$\mathbb{Z}_K[\textbf{\emph{v}},\textbf{\emph{v}}^{-1}]$ of the form
\begin{equation}\label{s_E}
s_E=\xi_E N_E \prod_{i \in I_E} \Psi_{E,i}(M_{E,i})
\end{equation}
where
\begin{enumerate}[(a)]
    \item $\xi_E$ is an element of $\mathbb{Z}_K$,
    \item $N_E= \prod_{\mathcal{C},j} v_{\mathcal{C},j}^{b_{\mathcal{C},j}}$ is a monomial in $\mathbb{Z}_K[\textbf{\emph{v}},\textbf{\emph{v}}^{-1}]$
          with $\sum_{j=0}^{e_\mathcal{C}-1}b_{\mathcal{C},j}=0$
          for all $\mathcal{C} \in \mathcal{A}/W$,
    \item $I_E$ is an index set,
    \item $(\Psi_{E,i})_{i \in I_E}$ is a family of $K$-cyclotomic polynomials in one variable,
    \item $(M_{E,i})_{i \in I_E}$ is a family of monomials in $\mathbb{Z}_K[\textbf{\emph{v}},\textbf{\emph{v}}^{-1}]$ such that
          if $M_{E,i} = \prod_{\mathcal{C},j} v_{\mathcal{C},j}^{a_{\mathcal{C},j}}$,
          then $\textrm{\emph{gcd}}(a_{\mathcal{C},j})=1$
          and $\sum_{j=0}^{e_\mathcal{C}-1}a_{\mathcal{C},j}=0$
          for all $\mathcal{C} \in \mathcal{A}/W$.
\end{enumerate}
\end{thm}

Equation (\ref{s_E}) gives the factorisation of $s_E$ into irreducible factors. 
The monomials $(M_{E,i})_{i \in I_E}$ are unique up to inversion, and we will call them \emph{potentially essential} for $W$.

\begin{rem}{\rm Theorem \ref{Schur element generic} was independently obtained by Rouquier \cite[Theorem 3.5]{Ro2} using a general argument on rational Cherednik algebras.}
\end{rem}

\begin{exmp}{\rm  Let us consider the example of $\mathfrak{S}_3$, which is isomorphic to $G(1,1,3)$.
We have
\begin{center}
$s_{E^{(3)}}=\Phi_2(v_{\textbf{t},0}v_{\textbf{t},1}^{-1})\Phi_3(v_{\textbf{t},0}v_{\textbf{t},1}^{-1}),\,s_{E^{(2,1)}}=v_{\textbf{t},0}^{-1}v_{\textbf{t},1}\Phi_3(v_{\textbf{t},0}v_{\textbf{t},1}^{-1}),\,s_{E^{(1,1,1)}}=v_{\textbf{t},0}^{-3}v_{\textbf{t},1}^{3}\Phi_2(v_{\textbf{t},0}v_{\textbf{t},1}^{-1})\Phi_3(v_{\textbf{t},0}v_{\textbf{t},1}^{-1}).$
\end{center}
}
\end{exmp}

 Let ${\varphi_{\textbf{m}}}: v_{\mathcal{C},j} \mapsto q^{m_{\mathcal{C},j}}$ be a cyclotomic specialisation. The canonical symmetrising trace on $\mathcal{H}(W)$ specialises via ${\varphi_{\textbf{m}}}$ to become the canonical symmetrising trace $\tau_{\varphi_{\textbf{m}}}$ on $\mathcal{H}_{\varphi_{\textbf{m}}}(W)$. The Schur elements of $\mathcal{H}_{\varphi_{\textbf{m}}}(W)$ with respect to $\tau_{\varphi_{\textbf{m}}}$ are $({\varphi_{\textbf{m}}}(s_E))_{E \in \Irr(W)}$, hence they can be written in the form (\ref{cycox}). We can again define functions $a^{\textbf{m}}: \Irr(W) \rightarrow \Z$ and $A^{\textbf{m}}: \Irr(W) \rightarrow \Z$ such that
$$a^{\textbf{m}}_E:=-{\rm val}_q({\varphi_{\textbf{m}}}(s_E)) \,\,\,\,\,\text{ and }\,\,\,\,\,A^{\textbf{m}}_E:=-{\rm deg}_q({\varphi_{\textbf{m}}}(s_E)).$$

\subsection{Families of characters and Rouquier families}\label{Roublo}

 Let ${\varphi_{\textbf{m}}}: v_{\mathcal{C},j} \mapsto q^{m_{\mathcal{C},j}}$ be a cyclotomic specialisation and let $\mathcal{H}_{\varphi_{\textbf{m}}}(W)$ be the corresponding cyclotomic Hecke algebra associated with $W$. How can we define families of characters for $\mathcal{H}_{\varphi_{\textbf{m}}}(W)$? We cannot apply Lusztig's original definition, because parabolic subgroups of complex reflection groups\footnote{The parabolic subgroups of a complex reflection group $W \subset {\rm GL}(\mathfrak{h})$ are the pointwise stabilisers of the subsets of $\mathfrak{h}$. It is a remarkable theorem by Steinberg \cite[Theorem 1.5]{St} that all parabolic subgroups of $W$ are again complex reflection groups.
}  do not have a nice presentation as in the real case, and certainly not a ``corresponding'' parabolic Hecke algebra. On the other hand,
we do not have a Kazhdan--Lusztig basis for $\mathcal{H}_{\varphi_{\textbf{m}}}(W)$, so we cannot construct Kazhdan--Lusztig cells and use them to define families of characters for complex reflection groups in the usual way. However, we can define the families of characters to be the Rouquier families of $\mathcal{H}_{\varphi_{\textbf{m}}}(W)$, that is, the blocks of $\mathcal{H}_{\varphi_{\textbf{m}}}(W)$ over the Rouquier ring $\mathcal{R}_K(q)$, where
 $$\mathcal{R}_K(q)=\Z_K[q,q^{-1},(q^n-1)^{-1}_{n \geq1}]. $$
Similarly to the real case, the Rouquier families are the non-empty subsets $B$ of $\Irr(W)$  that are minimal with respect to the property:
 $$\sum_{E \in B} \frac{\varphi_{\textbf{m}}(\chi_E(h))}{\varphi_{\textbf{m}}(s_E)} \in \mathcal{R}_K(q) \,\,\,\,\,\forall \, h \in \mathcal{H}(W).$$
 
 Brou\'e and Kim \cite{BK} determined the Rouquier families for the complex reflection groups of type $G(\ell,1,n)$, but their results are only true when $\ell$ is a power of a prime number or  ${\varphi_{\textbf{m}}}$ is a ``good'' cyclotomic specialisation. The same problem persists, and some new appear, in the determination of the Rouquier families for $G(\ell,p,n)$ by Kim \cite{Kim}. Malle and Rouquier \cite{MR} calculated the Rouquier families for some exceptional complex reflection groups and the dihedral groups, for a certain choice of cyclotomic specialisation. 
 More recently, we managed to determine the Rouquier families for all cyclotomic Hecke algebras of all complex reflection groups \cite{cras, berkeley, mybook, nagoya}, thanks to their property of ``semicontinuity'' (the term is due to C\'edric Bonnaf\'e). In order to explain this property, we will need some definitions.

Let $M = \prod_{\mathcal{C},j} v_{\mathcal{C},j}^{a_{\mathcal{C},j}}$ be a potentially essential monomial for $W$. We say that the family of integers $\textbf{m}=(m_{\mathcal{C},j})_{(\mathcal{C} \in
\mathcal{A}/W)(0\leq j \leq e_{\mathcal{C}}-1)}$ \emph{belongs to the potentially essential hyperplane} $H_M$ (of $\R^{\sum_{\mathcal{C}}e_\mathcal{C}}$)  if
$\sum_{\mathcal{C},j} m_{\mathcal{C},j}{a_{\mathcal{C},j}}=0$.

Suppose that $\textbf{m}$ belongs to no potentially essential hyperplane. Then the Rouquier families of $\mathcal{H}_{\varphi_{\textbf{m}}}(W)$ are called \emph{Rouquier families associated with no essential hyperplane}. 
Now suppose that  $\textbf{m}$ belongs to a unique potentially essential hyperplane $H$. Then the Rouquier families of $\mathcal{H}_{\varphi_{\textbf{m}}}(W)$ are called \emph{Rouquier families associated with} $H$. If they do not coincide with the  Rouquier families associated with no essential hyperplane, then $H$ is called an \emph{essential hyperplane} for $W$. All these notions are well-defined and they do not depend on the choice of $\textbf{m}$ because of the following theorem \cite[\S 4.4]{mybook}.

\begin{thm}{\rm  \textbf{(Semicontinuity  property of  Rouquier families)} } Let $\textbf{\emph{m}}=(m_{\mathcal{C},j})_{(\mathcal{C} \in \mathcal{A}/W)(0\leq j \leq e_{\mathcal{C}}-1)}$ be a family of integers and let  ${\varphi_{\textbf{\emph{m}}}}: v_{\mathcal{C},j} \mapsto q^{m_{\mathcal{C},j}}$ be the corresponding cyclotomic specialisation. The Rouquier families of $\mathcal{H}_{\varphi_{\textbf{\emph{m}}}}(W)$ are unions of the Rouquier families associated with the essential hyperplanes that $\textbf{\emph{m}}$ belongs to and they are minimal with respect to that property.
\end{thm}

Thanks to the above result, it is enough to do calculations in a finite number of cases in order to obtain the families of characters for all cyclotomic Hecke algebras, whose number is infinite.

\begin{exmp}{\rm  For $W=\mathfrak{S}_3$, the Rouquier families associated with no essential hyperplane are trivial. The hyperplane $H_M$ corresponding to the monomial $M=v_{\textbf{t},0}v_{\textbf{t},1}^{-1}$ is essential, and it is the unique essential hyperplane for $\mathfrak{S}_3$. Let  ${\varphi_{\textbf{m}}}: v_{\textbf{t},j} \mapsto q^{m_{j}}$, $j=0,1$, be a cyclotomic specialisation.
We have that $\textbf{m}=(m_0,m_1)$ belongs to $H_M$ if and only if $m_0=m_1$. There is a single Rouquier family associated with $H_M$, which contains all irreducible representations of 
$\mathfrak{S}_3$.}
\end{exmp}

We have also shown that the functions $a$ and $A$ are constant on the Rouquier families, for all cyclotomic Hecke algebras of all complex reflection groups \cite{degval, berkeley, nagoya}.

\subsection{Canonical basic sets} Given a cyclotomic Hecke algebra  $\mathcal{H}_{\varphi_{\textbf{m}}}(W)$ and a ring homomorphism $\theta: q \mapsto \eta \in \C\setminus\{0\}$, we obtain a semisimplicity criterion and  a decomposition map exactly as in \S\ref{canbasset}. A canonical basic set with respect to $\theta$ is also defined in the same way.

In \cite{chja1}, we showed the existence of canonical basic sets with respect to any $\theta$ for all cyclotomic Hecke algebras associated with finite Coxeter groups, 
that is, when the weight function $L$ in the definition of $\mathcal{H}(W,L)$ is also allowed to take negative values.

For non-real complex reflection groups, things become more complicated. For $W=G(\ell,1,n)$, consider the specialised Ariki--Koike algebra with relations 
\begin{equation}\label{ArikiKoike}
(\textbf{s} - \zeta_e^{s_0})(\textbf{s} - \zeta_e^{s_1})\cdots (\textbf{s}- \zeta_e^{s_{\ell - 1}}) = 0, \qquad 
(\textbf{t}_i-\zeta_e)(\textbf{t}_i+1) = 0\,\,\,\,\text{ for }i=1,\ldots,n-1.
\end{equation}
where $(s_0, \ldots , s_{\ell - 1}) \in \Z^{\ell}$ and $e \in \Z_{>0}$. With the use of Ariki's Theorem \cite{ar} and  Uglov's work on canonical bases for higher level Fock spaces \cite{ug}, Geck and Jacon \cite{GJ, jac1, jaca, GJlivre} have shown that, for a suitable choice of $\textbf{m}$, the corresponding function $a^{\textbf{m}}$ yields a canonical basic set for the above specialised Ariki--Koike algebra. This canonical basic set consists of the so-called ``Uglov $\ell$-partitions'' \cite[Definition 3.2]{jaca}. 
However, this does not work the other way round:  not all cyclotomic Ariki--Koike algebras admit canonical basic sets. For a study about which values of $\textbf{m}$ yield canonical basic sets, see \cite{Gerber}.

In \cite{chja2}, building on work by Genet and Jacon \cite{geja}, we generalised the above result to obtain canonical basic sets for all groups of type $G(\ell,p,n)$ with $n>2$, or $n=2$ and $p$ odd.

Finally, for the exceptional complex reflection groups of rank $2$ ($G_4$,\ldots,$G_{22}$), we have shown the existence of canonical basic sets for the cyclotomic Hecke algebras appearing in \cite{BM} with respect to any $\theta$ \cite{chmi}.

\section{Symplectic Reflection Algebras}

Let $V$ be a complex vector space of finite dimension $n$, and let  $G \subset  \mathrm{GL}(V)$ be a finite group. 
Let $\C[V]$ be the set of regular functions on $V$, which is the same thing as the symmetric algebra ${\rm Sym}(V^*)$ of the dual space of $V$.
The group $G$ acts on $\C[V]$ as follows:
$${}^g f(v):=f(g^{-1}v)\quad \forall  \,g\in G,\,f \in \C[V],\,v\in V.$$
We set
$$\C[V]^G:=\{f \in \C[V]\,|\,{}^g f=f \,\,\,\forall \,g \in G\},$$
the space of fixed points of $\C[V]$ under the action of $G$.
It is a classical problem in algebraic geometry to try and understand as a variety the space $$V/G = {\rm Spec}\,\C[V]^G.$$ 
Is the space $V/G$ singular? How much? The first question is answered by the following result, due to Shephard--Todd \cite{ShTo} and Chevalley \cite{Che}.

\begin{thm}\label{chevalley} The following statements are equivalent:
\begin{enumerate}[(1)]
\item $V/G$ is smooth. \smallbreak
\item $\C[V]^G$ is a polynomial algebra, on $n$ homogeneous generators. \smallbreak
\item $G$ is a complex reflection group.
\end{enumerate}
\end{thm}

\begin{exmp}{\rm Let $\mathfrak{S}_n$ act on $V=\C^n$ by permuting the coordinates.
Let $\C[V]=\C[X_1,\ldots,X_n]$ and let $\Sigma_1,\Sigma_2,\ldots,\Sigma_n$ be the elementary symmetric polynomials in $n$ variables.
We have
$$ \C[V]^{\mathfrak{S}_n}=\C[\Sigma_1(X_1,\ldots,X_n), \Sigma_2(X_1,\ldots,X_n),\ldots,\Sigma_n(X_1,\ldots,X_n)] .$$
More generally, we have
 $$ \C[V]^{G(\ell,1,n)}=\C[\Sigma_1(X_1^{\ell},\ldots,X_n^{\ell}), \Sigma_2(X_1^{\ell},\ldots,X_n^{\ell}),\ldots,\Sigma_n(X_1^{\ell},\ldots,X_n^{\ell})], $$
 where $G(\ell,1,n) \cong (\Z/ \ell \Z )^n \rtimes \mathfrak{S}_n$ and $(\Z/ \ell \Z )^n$  acts on $V$ by multiplying the coordinates by $\ell$-th roots of unity.
Note that $G(\ell,1,n)$ acts irreducibly on $V$ if and only if $\ell>1$.}
\end{exmp}

\begin{exmp}\label{examplemu2}{\rm Let $V=\C^2$ and let $G$ be a finite subgroup of ${\rm SL}_2(\C)$. Then $G$ is not a complex reflection group (in fact, it contains no pseudo-reflections at all).
The singular space $\C^2/G$ is called a \emph{Kleinian} (or \emph{Du Val}) singularity. 
The simplest example we can take is $$G=\left\{ \left(\begin{array}{cc} 1&0 \\0&1 \end{array}\right), \left(\begin{array}{rr} -1&0 \\0&-1 \end{array}\right)\right\} \cong \Z/2\Z \,;$$
we will use it to illustrate further notions.
}
\end{exmp}

\subsection{Symplectic reflection groups}\label{symreflgroups} The group $G$ in Example \ref{examplemu2} might not be a complex reflection group, but it is a symplectic reflection group, which is quite close. 
Moreover, the space  $\C^2/G$ is not smooth (following Theorem \ref{chevalley}), but it admits a symplectic resolution. 

Let  $(V,\omega_V)$ be a symplectic vector space, let  ${\rm Sp}(V)$ be the group of symplectic transformations on $V$ and let $G \subset {\rm Sp}(V)$ be a finite group.
The triple $(G,V,\omega_V)$ is called a \emph{symplectic triple}. A symplectic triple is {\em indecomposable}  if there is no $G$-equivariant splitting $V = V_1 \oplus V_2$ with $\omega_V(V_1, V_2) = 0$. Any  symplectic triple is a direct sum of indecomposable symplectic triples.

\begin{defn}\label{sympresn}{\rm Let $(G,V,\omega_V)$ be a symplectic triple and let $ (V/G)_{\rm sm}$ denote the smooth part of $V/G$.
A {\em symplectic resolution} of $V/G$ is a resolution of singularities $\pi : X \rightarrow V/G$ such that there exists a complex symplectic form $\omega_X$ on $X$ for which the isomorphism
$$\pi |_{\pi^{-1}((V/G)_{\rm sm})} : \pi^{-1}((V/G)_{\rm sm}) \rightarrow (V/G)_{\rm sm}$$
is a symplectic isomorphism.
}
\end{defn}

The existence of a symplectic resolution for $V/G$ is a very strong condition and implies that
the map $\pi$ has some very good properties, e.g., $\pi$ is ``semi-small'' \cite[Theorem 2.8]{ver}. Moreover, all crepant resolutions of $V/G$ are symplectic \cite[Theorem 2.5]{ver}.  

Verbitsky has shown that if $V/G$ admits a symplectic resolution, then $G$ is generated by symplectic reflections  \cite[Theorem 3.2]{ver}.

\begin{defn}\label{symplecticgroup}{\rm A \emph{symplectic reflection} is a non-trivial element $s \in {\rm Sp}(V)$ such that ${\rm rank}(s-$id$_V)=2$.
The symplectic triple $(G,V,\omega_V)$ is a \emph{symplectic reflection group} if $G$  is generated by  symplectic reflections.}
\end{defn}

Hence, if the space $V/G$ admits a symplectic resolution, then $(G,V,\omega_V)$ is a symplectic reflection group; the converse is not true. The classification of such symplectic reflection groups is almost complete thanks to the representation theory of symplectic reflection algebras.

\begin{exmp}{\rm Following Example \ref{examplemu2},
 let $G$ be the cyclic group of order $2$, denoted by $\mu_2$, acting on $V = \C \oplus \C^*$ by multiplication by $-1$. Let $\omega_V$ be the \emph{standard symplectic form} on $V$, that is,
 \begin{equation}\label{standardsymform}
\omega_V(y_1 \oplus x_1, y_2 \oplus x_2) = x_2(y_1)-x_1(y_2).
\end{equation} 
  Letting $\C[V ] = \C[X ,Y]$, we see that $\C[V ]^G = \C[X^2, XY, Y^2] \cong \C[A, B, C]/(AC - B^2)$, the quadratic cone. This has an isolated singularity at the origin, \emph{i.e.,} at the zero orbit, which can be resolved by blowing up there. The resulting resolution $\pi : T^*\mathbb{P}^1  \rightarrow V/G$ is a symplectic resolution where $T^*\mathbb{P}^1$ has its canonical symplectic structure. }
  \end{exmp}

The classification of  (indecomposable) symplectic reflection groups is due to Huffman--Wales \cite{95}, Cohen \cite{37}, and Guralnick--Saxl \cite{87}. Except for a finite list of explicit exceptions with ${\rm dim}_{\C}(V)\leq 10$, there are two classes of symplectic reflection groups:
\begin{itemize}
\item {\bf Wreath products.} Let $\Gamma \subset {\rm SL}_2 (\C)$ be finite: such groups are called Kleinian subgroups and they preserve the canonical symplectic structure on $\C^2$. Set
$$V =\underbrace{\C^2 \oplus \C^2 \oplus \cdots \oplus\C^2}_{n \text{ summands}}$$
with the symplectic form $\omega_V$ induced from that on $\C^2$ and let  $G = \Gamma \wr \mathfrak{S}_n$ act in the obvious
way on $V$.\smallbreak
\item {\bf Complex reflection groups.} Let $G \subset {\rm GL} (\mathfrak{h})$ be a complex reflection group. Set $V = \mathfrak{h} \oplus \mathfrak{h}^*$  with its standard symplectic form $\omega_V$ (see (\ref{standardsymform})) and with $G$ acting diagonally. 
\end{itemize}
In both of the above cases,  $(G,V,\omega_V)$ is a symplectic reflection group. 

\begin{rem}{\rm Note that in the second case,  where $G$ is a complex reflection group, the space $\mathfrak{h}/G$ is smooth, but $V/G$ is not. The symplectic reflections in $(G,V,\omega_V)$ are the pseudo-reflections in $(G,\mathfrak{h})$. 
}
\end{rem}

\begin{rem}{\rm There is a small overlap between the two main families of symplectic reflection groups, namely the complex reflection groups of type $G(\ell,1,n)$.}
\end{rem}

In \cite[\S 1.3 and  \S 1.4]{Wa} Wang observes that if  $G = \Gamma \wr \mathfrak{S}_n$ for some $\Gamma \subset {\rm SL}_2 (\C)$, then $V/G$ has a symplectic resolution given by the Hilbert scheme of $n$ points on the minimal resolution of the Kleinian singularity $\C^2/\Gamma$. In Section \ref{lastsection} we will see what happens in the case where $G $ is a complex reflection group.

\subsection{The symplectic reflection algebra ${\bf H}_{t,\bc}(G)$} From now on, let $(G,V,\omega_V)$ be a symplectic reflection group and let $\mathcal{S}$ be the set of all symplectic reflections in $G$.

\begin{defn}\label{skew}{\rm The \emph{skew-group ring} $\C[V] \rtimes G$ is, as a vector space, equal to $\C[V] \otimes \C G$ and the multiplication is given by  
$$ g \cdot f = {}^g f \cdot g \,\,\,\,\,\forall\,g \in G,\,f \in \C[V]. $$
}
\end{defn}

The centre $Z(\C[V] \rtimes G)$ of the skew-group ring is equal to $\C[V]^G$. It has  been an insight of Etingof and Ginzburg \cite{EG}, which goes back to (at least) Crawley-Bovey and Holland \cite{CBH}, that, in order to understand ${\rm Spec}\,\C[V]^G$, we could look at deformations of $\C[V] \rtimes G$, hoping that the centre of the deformed algebra is itself a deformation of $\C[V]^G$. These deformations are the symplectic reflection algebras.

Let $s \in \mathcal{S}$. The spaces ${\rm Im}(s - {\rm id}_V)$ and ${\rm Ker}(s - {\rm id}_V)$ are symplectic subspaces of $V$ with ${\rm dim}_{\C}{\rm Im}(s - {\rm id}_V)= 2$ and $V = {\rm Im}(s - {\rm id}_V) \oplus  {\rm Ker}(s - {\rm id}_V)$. Let $\omega_s$ be the $2$-form on $V$ whose
 restriction to  ${\rm Im}(s - {\rm id}_V)$ is $\omega_V$ and whose restriction to  ${\rm Ker}(s - {\rm id}_V)$ is zero. 
Let $\omega_{V^*}$ be the symplectic form on $V^*$  corresponding to $\omega_V$ (under the identification of $V$ and $V^*$
induced by $\omega_V$), and let $TV^{*}$ denote the tensor algebra on $V^*$. Finally, let $\bc: \mathcal{S} \rightarrow \C$ be a {\em conjugacy invariant function}, that is, a map such that
$$\bc(gsg^{-1})=\bc(s) \,\,\,\,\,\forall \, s \in \mathcal{S},\,g\in G.$$
 
\begin{defn}\label{syreal}
 {\rm  Let $t \in \C$. We define the \emph{symplectic reflection algebra} ${\bf H}_{t,\bc}(G)$ of $G$ to be
$$
{\bf H}_{t,\bc}(G) := TV^{*} \rtimes G/ \langle [u, v] - (t\,\omega_{V^*}(u, v) - 2 \sum_{s \in \mathcal{S}} \bc(s)\,\omega_s(u, v)\, s ) \,|\, u, v \in V^* \rangle.
 $$}
\end{defn}

Note that the above definition simply describes how two vectors in $V^*$ commute with each other in ${\bf H}_{t,\bc}(G)$, and that we have $[u,v] \in \C G$ for all $u,v \in V^*$.

\begin{rem}{\rm For all $\lambda \in \C^\times$, we have ${\bf H}_{\lambda t,\lambda\bc}(G) \cong {\bf H}_{t,\bc}(G)$. So we only need to consider the cases $t=1$ and $t=0$.
}
\end{rem}

\begin{rem}{\rm We have ${\bf H}_{0,0}(G) = \C[V] \rtimes G$.
}
\end{rem}

\begin{exmp}\label{18}
{\rm Let us consider the example of the cyclic group $\mu_2=\langle s \rangle$ acting on $V=\C^2$, so that 
\begin{center}
$sx=-x$, \,$sy=-y$\, and \,$\omega_{V^*}(y,x)$=1,
\end{center}
where $\{x,y\}$ is a basis of $(\C^2)^*$. We have $\omega_s = \omega_{V^*}$, since ${\rm Im}(s - {\rm id}_V)=V$.  Then ${\bf H}_{t,\bc}(\mu_2)$ is the quotient of $ \C\langle x,y,s \rangle$ by the relations:
$$s^2=1,\,\,sx=-xs,\,\,sy=-ys,\,\,[y,x]=t-2\bc(s) s.$$ 
}
\end{exmp}

\begin{exmp}{\rm Let $V=\C^2$. Then ${\rm Sp}(V)={\rm SL}_2(\C)$ and we can take $G$ to be any finite subgroup of ${\rm SL}_2(\C)$.  Let $\{x,y\}$ be a basis of $(\C^2)^*$ such that $\omega_{V^*}(y, x) = 1$. Every $g \neq 1$ in $G$ is a symplectic reflection and $\omega_g = \omega_{V^*}$. Then
$${\bf H}_{t,\bc}(G) =\C\langle x,y \rangle \rtimes G / \langle [y,x] -(t-2 \sum_{g \in G\setminus\{1\}}\bc(g)g) \rangle.$$
}
\end{exmp}

There is a natural filtration $\mathcal{F}$ on ${\bf H}_{t,\bc}(G)$ given by putting $V^*$ in degree one and $G$ in degree zero. The crucial result by Etingof and Ginzburg is the Poincar\'e--Birkhoff--Witt (PBW) Theorem \cite[Theorem 1.3]{EG}. 

\begin{thm}  There is an isomorphism of algebras 
$${\rm gr}_{\mathcal{F}} ({\bf H}_{t,\bc}(G) ) \cong \C[V ] \rtimes G,$$
given by $\sigma(v) \mapsto v$, $\sigma(g) \mapsto g$, where $\sigma(h)$ denotes the image of $h \in {\bf H}_{t,\bc}(G)$ in ${\rm gr}_{\mathcal{F}} ({\bf H}_{t,\bc}(G) )$. In particular, there is an isomorphism of vector spaces $${\bf H}_{t,\bc}(G) \cong \C[V ] \otimes \C G.$$ Moreover, symplectic reflection algebras are the only deformations of $\C[V ] \rtimes G$
with this property (PBW property).
\end{thm}

The most important consequence of the PBW Theorem is that it gives us an explicit basis of the symplectic reflection algebra. The proof of it is an application of a general result by Braverman and Gaitsgory: If $I$ is a two-sided ideal of $TV^* \rtimes G$ generated by a space $U$ of elements of degree at most two, then \cite[Theorem 0.5]{BrGa} gives necessary and sufficient conditions so that the quotient $TV^* \rtimes G/I$ has the PBW property. The PBW property also implies that ${\bf H}_{t,\bc}(G)$ has some good ring-theoretic properties, for example:

\begin{cor} 
\begin{enumerate}[(i)] 
\item  The algebra ${\bf H}_{t,\bc}(G)$  is a Noetherian ring. \smallbreak
\item  ${\bf H}_{t,\bc}(G)$ has finite global dimension.
\end{enumerate}
\end{cor}

\begin{rem}{\rm
For general pairs $(G, V )$ a description of PBW deformations of $\C[V ] \rtimes G$ was originally given by Drinfeld \cite{44}. In the symplectic case this was rediscovered by Etingof and Ginzburg as above, and Drinfeld's general case was described in detail by Ram and Shepler \cite{124}. 
}\end{rem}

\subsection{The spherical subalgebra}
We saw in the previous subsection that the skew-group ring $\C[V] \rtimes G$ is not commutative and that its centre $Z(\C[V] \rtimes G)$ is equal to $\C[V]^G$. We will now see that $\C[V] \rtimes G$ contains another subalgebra isomorphic to $\C[V]^G$.

Let $e:=\frac{1}{|G|} \sum_{g\in G}g$
be the trivial idempotent in $\C G$. One can easily check that the map
\begin{equation}\label{cvg}
\begin{array}{ccc}
\C[V]^G & \rightarrow & e (\C[V] \rtimes G) e \\
f & \mapsto &e f e
\end{array}
\end{equation}
is an algebra isomorphism. We have $efe=fe$, for all $f \in \C[V]^G$.

\begin{defn} {\rm We define the \emph{spherical subalgebra} of ${\bf H}_{t,\bc}(G)$ to be the algebra
$${\bf U}_{t,\bc}(G):=e\, {\bf H}_{t,\bc}(G)\, e.$$}
\end{defn}

The filtration $\mathcal{F}$ on ${\bf H}_{t,\bc}(G)$ induces, by restriction, a filtration on ${\bf U}_{t,\bc}(G)$.
The PBW Theorem, in combination with (\ref{cvg}), implies that there is an isomorphism of algebras
$${\rm gr}_{\mathcal{F}} ({\bf U}_{t,\bc}(G) )\cong e (\C[V] \rtimes G) e \cong  \C[V ] ^G  $$
and an isomorphism of vector spaces
$$ {\bf U}_{t,\bc}(G) \cong  \C[V ] ^G .$$
Thus, the spherical subalgebra provides a  flat deformation of the coordinate ring of $V/G$, as desired.

\begin{exmp}\label{likelie}{\rm 
Let $G = \mu_2=\langle s \rangle$ acting on $V= \C^2$ as in Example \ref{18}. Then $e=\frac{1}{2}(1+s)$. The spherical subalgebra ${\bf U}_{t,\bc}(\mu_2)$ is generated as a $\C$-algebra by
$${\bf h}: = -\frac{1}{2} e(xy + yx)e,\,\, {\bf e}: = \frac{1}{2}  ex^2e\,\,\text{ and }\,\, {\bf f}: = \frac{1}{2} ey^2e.$$ There are relations
$$[{\bf e}, {\bf f} ] = t{\bf h},\,\, [{\bf h}, {\bf e}] = -2t{\bf e},\,\, [{\bf h}, {\bf f} ] = 2t{\bf f} \,\,\text{ and }\,\, {\bf e}{\bf f} = (2\bc(s) - {\bf h}/2)(t/2 - \bc(s) - {\bf h}/2).$$
So if $t = 0$,  ${\bf U}_{t,\bc}(\mu_2)$ is commutative, while if $t = 1$,  ${\bf U}_{t,\bc}(\mu_2)$ is a central quotient of the
enveloping algebra of $\mathfrak{sl}_2(\C)$.
}
\end{exmp}

The space ${\bf H}_{t,\bc}(G)e$ is a $({\bf H}_{t,\bc}(G),{\bf U}_{t,\bc}(G))$-bimodule and it is called {\em the Etingof--Ginzburg sheaf}.
The following result is known as the ``double centraliser property'' \cite[Theorem 1.5]{EG}.

\begin{prop} \begin{enumerate}[(i)]
\item The right ${\bf U}_{t,\bc}(G)$-module ${\bf H}_{t,\bc}(G)e$ is reflexive. \smallbreak
\item ${\rm End}_{{\bf H}_{t,\bc}(G)}({\bf H}_{t,\bc}(G)e)^{op} \cong {\bf U}_{t,\bc}(G)$. \smallbreak
\item ${\rm End}_{{\bf U}_{t,\bc}(G)^{op}} ({\bf H}_{t,\bc}(G)e)  \cong {\bf H}_{t,\bc}(G)$.
\end{enumerate}
\end{prop}

This is important, because, in general, we have an explicit presentation of ${\bf H}_{t,\bc}(G)$, but not of ${\bf U}_{t,\bc}(G)$. The above result allows us to study ${\bf U}_{t,\bc}(G)$ by studying ${\bf H}_{t,\bc}(G)$ instead.  

\subsection{The centre of ${\bf H}_{t,\bc}(G)$}
The behaviour of the centre of the spherical subalgebra observed in Example \ref{likelie} is the same for all symplectic reflection groups \cite[Theorem 1.6]{EG}.

\begin{thm} \begin{enumerate}[(i)]
\item If $t=0$, then ${\bf U}_{t,\bc}(G)$ is commutative. \smallbreak
\item If $t \neq 0$, then $Z({\bf U}_{t,\bc}(G))=\C$.
\end{enumerate}
\end{thm}

Now the double centraliser property can be used to prove the following result relating the centres of ${\bf U}_{t,\bc}(G)$ and ${\bf H}_{t,\bc}(G)$.
\begin{thm}{\rm  \textbf{(The Satake isomorphism)} } The map $z \mapsto z  e$ defines an algebra isomorphism
$Z({\bf H}_{t,\bc}(G)) \cong Z({\bf U}_{t,\bc}(G))$ for all parameters $(t,\bc)$.
\end{thm}

\begin{cor} \begin{enumerate}[(i)]
\item If $t=0$, then $Z({\bf H}_{t,\bc}(G)) \cong {\bf U}_{t,\bc}(G)$. \smallbreak
\item If $t \neq 0$, then $Z({\bf H}_{t,\bc}(G))=\C$.
\end{enumerate}
\end{cor}

Thus, the symplectic reflection algebra  ${\bf H}_{t,\bc}(G)$ produces a commutative deformation of the space $V/G$ when $t = 0$.

\subsection{Symplectic resolutions} In this subsection, we will focus on the case $t=0$. Set $Z_\bc(G):=Z({\bf H}_{0,\bc}(G))$.
We have  $Z_\bc(G) \cong {\bf U}_{0,\bc}(G)$, and so ${\bf H}_{0,\bc}(G)$ is a finitely generated $Z_\bc(G)$-module.

\begin{defn}\label{CalogeroMoser} {\rm The 
\emph{generalised Calogero--Moser space} $X_{\bc}(G)$ is defined to be the affine variety
${\rm Spec}\,Z_\bc(G)$.}
\end{defn}

Since the associated
graded of $Z_\bc(G)$ is $\C[V ]^G$ (with respect to the filtration $\mathcal{F}$), $X_{\bc}(G)$ is irreducible.
The following result, due to Ginzburg--Kaledin \cite[Proposition 1.18 and Theorem 1.20]{gika} and Namikawa \cite[Corollary 2.10]{na}, gives us a criterion for $V/G$ to admit a symplectic resolution, using the geometry of the generalised Calogero--Moser space.

\begin{thm}\label{smoothness} Let $(G,V, \omega_V)$ be an (irreducible) symplectic reflection group. 
The space $V/G$ admits a symplectic resolution if and only if $X_{\bc}(G)$ is smooth for generic values of $\bc$ (equivalently, there exists $\bc$ such that  $X_{\bc}(G)$ is smooth).
\end{thm}

\begin{exmp}{\rm
Consider again the example of $\mu_2=\langle s \rangle$ acting on $\C^2$. The centre of ${\bf H}_{0,\bc}(\mu_2)$ is generated by $A := x^2$, $B := xy - \bc(s)s$
and $C: = y^2$. Thus, $$X_{\bc}(\mu_2) \cong \C[A,B,C]/(AC - (B+\bc(s))(B-\bc(s)))$$
is the affine cone over $\mathbb{P}^1 \subset \mathbb{P}^2$ when $\bc(s)=0$, but is a smooth affine surface for  $\bc(s)\neq 0$. }
\end{exmp}

As we mentioned in Subsection \ref{symreflgroups}, if $G= \Gamma \wr \mathfrak{S}_n$ for some $\Gamma \subset {\rm SL}_2 (\C)$, then $V/G$ always admits a symplectic resolution, that is, $X_{\bc}(G)$ is smooth for generic $\bc$. On the other hand, if $G \subset {\rm GL}(\mathfrak{h})$ is a complex reflection group acting on $V = \mathfrak{h} \oplus \mathfrak{h}^*$, this is not always the case. Etingof and Ginzburg proved that  $X_{\bc}(G)$ is smooth for generic $\bc$ when $G=G(\ell,1,n)$
 \cite[Corollary 1.14]{EG}. However, Gordon showed that, for most finite Coxeter groups not of type $A_n$ or $B_n$, $X_{\bc}(G)$ is a singular variety for all choices of the parameter $\bc$ \cite[Proposition 7.3]{go}. Finally, using the Calogero--Moser partition of $\Irr(G)$ described in \cite{GM}, Bellamy proved that $X_{\bc}(G)$ is smooth for generic values of $\bc$ if and only if $G=G(\ell,1,n)$ or $G=G_4$  \cite[Theorem 1.1]{bel}. We will revisit this result in Section \ref{lastsection}.

Following the classification of symplectic reflection groups, and all the works mentioned above, the classification of quotient singularities admitting symplectic resolutions is (almost) complete.

\subsection{Rational Cherednik algebras} From now on, let $W  \subset {\rm GL}(\mathfrak{h})$ be a complex reflection group and let $V=\mathfrak{h} \oplus \mathfrak{h}^*$. There is a natural pairing $(\,,\,):\mathfrak{h} \times \mathfrak{h}^* \rightarrow \C$ given by $(y,x):=x(y)$. Then the standard symplectic form $\omega_V$ on $V$ is given by
$$\omega_V(y_1\oplus x_1,y_2 \oplus x_2)=(y_1,x_2)-(y_2,x_1).$$
The  triple $(W,V,\omega_V)$ is a symplectic reflection group. The set $\mathcal{S}$ of all symplectic reflections in $(W,V,\omega_V)$ coincides with the set of pseudo-reflections in $(W,\mathfrak{h})$. Let $\bc : \mathcal{S} \rightarrow \C$ be a conjugacy invariant function.

\begin{defn}
{\rm The \emph{rational Cherednik algebra} of $W$ is the symplectic reflection algebra ${\bf H}_{t,\bc}(W)$ associated to $(W,V,\omega_V)$.}
\end{defn}

For $s \in \mathcal{S}$, fix $\alpha_s \in  \mathfrak{h}^*$ to be a basis of the 
one-dimensional vector space ${\rm Im}(s - {\rm id}_V)|_{\mathfrak{h}^*}$ and $\alpha_s^\vee \in  \mathfrak{h}$ to be a basis of the 
one-dimensional vector space ${\rm Im}(s - {\rm id}_V)|_{\mathfrak{h}}$. 
Then ${\bf H}_{t,\bc}(W)$ is the quotient of $TV^* \rtimes W$ by the relations:
\begin{equation}\label{cherednikrel}
[x_1,x_2]=0,\,\,[y_1,y_2]=0,\,\,[y,x]=t (y,x) - 2 \sum_{s \in \mathcal{S}} \bc(s) \,\frac{(y,\alpha_s)(\alpha_s^\vee,x)}{(\alpha_s^\vee,\alpha_s)}\,s
\end{equation}
for all $x_1,x_2,x \in \mathfrak{h}^*$ and $y_1,y_2,y \in \mathfrak{h}$.

\begin{exmp}{\rm Let $W=\mathfrak{S}_n$ and $\mathfrak{h}=\C^n$. 
Choose a basis $x_1,\ldots,x_n$ of $\mathfrak{h}^*$ and a dual basis $y_1,\ldots,y_n$ of $\mathfrak{h}$ so that
$$\sigma x_i = x_{\sigma(i)}\sigma\,\,\,\text{ and }\,\,\, \sigma(y_i) = y_{\sigma(i)}\sigma \,\,\,\,\, \forall \,\sigma \in  \mathfrak{S}_n,\,1 \leq i \leq n.$$
The set $\mathcal{S}$ is the set of all transpositions in $\mathfrak{S}_n$. We denote by $s_{ij}$ the transposition $(i,j)$.
Set
$$\alpha_{ij}:=x_i-x_j \,\,\,\text{ and }\,\,\, \alpha_{ij}^\vee = y_i-y_j \,\,\,\,\, \forall \,1 \leq i < j \leq n.$$
We have $(\alpha_{ij}^\vee,\alpha_{ij})=2$. 
There is a single conjugacy class in $\mathcal{S}$, so take $\bc \in \C$. Then ${\bf H}_{t,\bc}(\mathfrak{S}_n)$ is the quotient of $TV^* \rtimes \mathfrak{S}_n$ by the relations:
$$[x_i,x_j]=0,\,\,\,[y_i,y_j]=0,\,\,\,[y_i,x_i]= t - \bc \sum_{j \neq i}s_{ij},\,\,\,[y_i,x_j]=\bc \,s_{ij} \,\,\,\text{for } i \neq j.$$
}
\end{exmp}

\section{Rational Cherednik Algebras at $t =1$} 

The PBW Theorem implies that the rational Cherednik algebra ${\bf H}_{1,\bc}(W)$, as a vector space,  has a ``triangular decomposition''
$${\bf H}_{1,\bc}(W) \cong \C[\mathfrak{h}] \otimes \C W \otimes \C[\mathfrak{h}^*].$$
Another famous example of a triangular decomposition is the one of the enveloping algebra $U(\mathfrak{g})$ of a finite dimensional, semisimple complex Lie algebra $\mathfrak{g}$  (into the enveloping algebras of the Cartan subalgebra, the nilpotent radical of the Borel subalgebra and its opposite). In the representation theory of $\mathfrak{g}$, one of the categories of modules most studied, and best understood, is category $\mathcal{O}$, the abelian category generated by all highest weight modules. Therefore, it makes sense to want to construct and study an analogue of category $\mathcal{O}$ for rational Cherednik algebras.  

\subsection{Category $\mathcal{O}$}\label{O} Let ${\bf H}_{1,\bc}(W)$-mod be the category of all finitely generated ${\bf H}_{1,\bc}(W)$-modules. We say that a module $M \in {\bf H}_{1,\bc}(W)$-mod is \emph{locally nilpotent} for the action of $\mathfrak{h} \subset \C[\mathfrak{h}^*]$ if for each $m \in M$ there exists $N > > 0$ such that $\mathfrak{h}^N \cdot m = 0$.

\begin{defn}
We define $\mathcal{O}$ to be the category of all finitely generated ${\bf H}_{1,\bc}(W)$-modules that are locally nilpotent for the action of  $\mathfrak{h} \subset \C[\mathfrak{h}^*]$.
\end{defn}

\begin{rem}{\rm Each module in category $\mathcal{O}$ is finitely generated as a $ \C[\mathfrak{h}]$-module.}
\end{rem}

Category $\mathcal{O}$ has been thoroughly studied in \cite{GGOR}. Proofs of all its properties presented here can be found in this paper.

For all $E \in  \Irr(W)$, we set
$$\Delta(E):={\bf H}_{1,\bc}(W) \otimes_{\C[\mathfrak{h}^*] \rtimes W}  E ,$$
where $\C[\mathfrak{h}^*]$ acts trivially on  $E$ (that is, the augmentation ideal $\C[\mathfrak{h}^*]_+$ acts on $E$ as zero) and $W$ acts naturally. The module $\Delta(E)$ belongs to $\mathcal{O}$ and is called a \emph{standard module} (or \emph{Verma module}). Each standard module $\Delta(E)$ has a simple head ${\rm L}(E)$ and the set 
$$\{ {\rm L}(E) \,| \, E \in \Irr(W)\}$$
is a complete set of pairwise non-isomorphic simple modules of the category $\mathcal{O}$. Every module in $\mathcal{O}$ has finite length, so we obtain a well-defined square decomposition matrix
$${\bf D}=([\Delta(E):{\rm L}(E')])_{E,E' \in \Irr(W)},$$
where $[\Delta(E):{\rm L}(E')]$ equals the multiplicity with which the simple module ${\rm L}(E')$ appears in the composition series of $\Delta(E)$. We have $[\Delta(E):{\rm L}(E)]=1$.

\begin{prop} The following are equivalent:
\begin{enumerate}[(1)]
\item $\mathcal{O}$ is semisimple. \smallbreak
\item $\Delta(E)={\rm L}(E)$ for all $E \in \Irr(W)$. \smallbreak
\item ${\bf D}$ is the identity matrix.
\end{enumerate}
\end{prop}

Now, there exist several orderings on the set of standard modules of $\mathcal{O}$ (and consequently on $\Irr(W)$) for which 
$\mathcal{O}$ is a highest weight category in the sense of \cite{CPS} (see also \cite[\S 5.1]{Ro2}). If  $<_{\mathcal{O}}$ is such an ordering on $\Irr(W)$, and if $[\Delta(E):{\rm L}(E')] \neq 0$ for some 
$E,E' \in \Irr(W)$, then either $E=E'$ or $E' <_{\mathcal{O}} E$. Thus, we can arrange the rows of ${\bf D}$ so that the decomposition matrix  is lower unitriangular.
We will refer to these orderings on $\Irr(W)$ as \emph{orderings on the category} $\mathcal{O}$.
A famous example of such an ordering  is the one given by the $c$-function.

\subsection{A change of parameters and the $c$-function}

In order to relate rational Cherednik algebras with cyclotomic Hecke algebras via the $\KZ$-functor in the next subsection, we need to change the parametrisation of $ {\bf H}_{1,\bc}(W)$. As in \S\ref{Hecke}, let $\mathcal{A}$ denote the set of reflecting hyperplanes of $W$. For $H \in \mathcal{A}$, let $W_H$ be the pointwise stabiliser of $H$ in $W$. The group $W_H$ is cyclic and its order, denoted by $e_{\mathcal{C}}$, only depends on the orbit $\mathcal{C} \in \mathcal{A}/W$ that $H$ belongs to. We have that
$$\mathcal{S} = \bigcup_{H \in \mathcal{A}} W_H\setminus \{1\}.$$
For each $s \in W_H\setminus \{1\}$, we have ${\rm Ker} \alpha_s=H$.
Without loss of generality, we may assume that $\alpha_s=\alpha_{s'}$ and $\alpha_s^\vee=\alpha_{s'}^\vee$ for all $s,\,s' \in W_H\setminus \{1\}$.
Set $\alpha_H:=\alpha_s$ and $\alpha_H^\vee:=\alpha_s^\vee$. Then the third relation in (\ref{cherednikrel}) becomes
$$[y,x]=(y,x) - 2 \sum_{H \in \mathcal{A}} \,\frac{(y,\alpha_H)(\alpha_H^\vee,x)}{(\alpha_H^\vee,\alpha_H)}\sum_{s \in W_H \setminus\{1\}} \bc(s)\,s \,\,\,\,\,\,\,\forall\, x \in \mathfrak{h}^*,\,y \in \mathfrak{h}.$$
We define a family of complex numbers ${\bf k}=(k_{\mathcal{C},j})_{(\mathcal{C} \in \mathcal{A}/W)(0 \leq j \leq e_{\mathcal{C}}-1)}$ by
$$-2 \sum_{s \in W_H \setminus\{1\}} \bc(s)\,s =  \sum_{s \in W_H\setminus\{1\}} \left(\sum_{j=0}^{e_\mathcal{C}-1} \det(s)^{-j} ({k}_{\mathcal{C},j} - {k}_{\mathcal{C},j-1})\right) s
\,\,\,\,\,\,\,\text{ for } H \in \mathcal{C}$$
with $k_{\mathcal{C},-1}=0$. This implies that
$$ \bc(s) = -\frac{1}{2} \sum_{j=0}^{e_\mathcal{C}-1} \det(s)^{-j} ({k}_{\mathcal{C},j} - {k}_{\mathcal{C},j-1}).$$ From now on, we will denote by
${\bf H}_{{\bf k}}(W)$ the quotient of $TV^* \rtimes W$ by the relations:
$$
[x_1,x_2]=0,\,\,[y_1,y_2]=0,\,\,[y,x]= (y,x) + \sum_{H \in \mathcal{A}} \,\frac{(y,\alpha_H)(\alpha_H^\vee,x)}{(\alpha_H^\vee,\alpha_H)} \,\gamma_H$$
where
$$\gamma_H= \sum_{w \in W_H\setminus\{1\}} \left(\sum_{j=0}^{e_\mathcal{C}-1} \det(w)^{-j} ({k}_{\mathcal{C},j} - {k}_{\mathcal{C},j-1})\right) w
$$
for all $x_1,x_2,x \in \mathfrak{h}^*$ and $y_1,y_2,y \in \mathfrak{h}$. We have ${\bf H}_{{\bf k}}(W)={\bf H}_{1,\bc}(W)$. \smallbreak

Let $E \in \Irr(W)$. We denote by $c_E$ the scalar  by which the element $$-\sum_{H \in \mathcal{A}} \sum_{j=0}^{e_{\mathcal{C}}-1} \left( \sum_{w \in W_H} (\det w)^{-j}w \right) k_{\mathcal{C},j}\in Z(\C W)$$ acts on $E$. We obtain thus a function $c: \Irr(W) \rightarrow \C,\, E \mapsto c_E$. The $c$-function defines an ordering $<_c$ on the category $\mathcal{O}$ as follows: For all $E,E' \in \Irr(W)$,
$$ E' <_c E \,\,\,\text{ if and only if }\,\,\, c_E - c_{E'} \in \Z_{>0}.$$

\begin{rem}{\rm If $c_E - c_{E'} \notin \Z \setminus\{0\}$ for all $E,E' \in \Irr(W)$, then ${\bf D}$ is the identity matrix, and thus $\mathcal{O}$ is semisimple.
}
\end{rem}

\begin{rem}{\rm In the rational Cherednik algebra literature the function $c$ is usually taken to be the negative of the one defined here. In the context of this paper the above definition is more natural. In both cases we obtain an ordering on the category $\mathcal{O}$.
}
\end{rem}

\subsection{The KZ-functor}\label{KZ} Following \cite[5.3]{GGOR}, there exists an exact factor, known as the \emph{Knizhnik--Zamalodchikov functor} or simply $\KZ$, between the category $\mathcal{O}$ of  ${\bf H}_{{\bf k}}(W)$ and the category of representations of a certain specialised Hecke algebra $\mathcal{H}_{{\bf k}}(W)$.  Using the notation of \S\ref{Hecke}, the specialised Hecke algebra $\mathcal{H}_{{\bf k}}(W)$ is a quotient of the group
algebra $\mathbb{C}B_W$ by the ideal
generated by the elements of the form
$$(\textbf{s}-{\rm exp}(2\pi i k_{\mathcal{C},0}))(\textbf{s}-\zeta_{e_\mathcal{C}}{\rm exp}(2\pi i k_{\mathcal{C},1})) \cdots  (\textbf{s}-\zeta_{e_\mathcal{C}}^{e_\mathcal{C}-1}{\rm exp}(2\pi i k_{\mathcal{C},e_\mathcal{C}-1})),$$
where $\mathcal{C}$ runs over the set $\mathcal{A}/W$ and
$\textbf{s}$ runs over the set of monodromy generators around 
the images in $\mathfrak{h}^{\textrm{reg}}/W$ of 
the elements of  $\mathcal{C}$. The algebra $\mathcal{H}_{{\bf k}}(W)$ is obtained from the generic Hecke algebra  $\mathbb{C}[{\bf v}, {\bf v}^{-1}]\mathcal{H}(W)$ 
via the specialisation $\Theta: v_{\mathcal{C},j}^{N_W} \mapsto {\rm exp}(2\pi i k_{\mathcal{C},j})$ (recall that $N_W$ is the power to which the indeterminates $ v_{\mathcal{C},j}$ appear in the defining relations of the generic Hecke algebra so that the algebra $\C({\bf v})\mathcal{H}(W)$ is split; see (\ref{NW})). 
We always assume that Hypothesis \ref{assumptions} holds for $\mathcal{H}(W)$.

The functor $\KZ$ is represented by a projective object  $P_{\KZ} \in \mathcal{O}$, and we have $\mathcal{H}_{{\bf k}}(W) \cong {\rm End}_{{\bf H}_{{\bf k}}(W)}(P_{\KZ})^{op}$
\cite[5.4]{GGOR}.
Based on this, we have the following result due to Vale \cite[Theorem 2.1]{Vale}:

\begin{prop} The following are equivalent:
\begin{enumerate}[(1)]
\item ${\bf H}_{{\bf k}}(W)$ is a simple ring. \smallbreak
\item $\mathcal{O}$ is semisimple.\smallbreak
\item$\mathcal{H}_{{\bf k}}(W)$ is semisimple. 
\end{enumerate}
\end{prop}

We can thus use the semisimplicity criterion  for $\mathcal{H}_{{\bf k}}(W)$
 given by Theorem \ref{sscriterion} in order to determine for which values of ${\bf k}$ the category $\mathcal{O}$ is semisimple.

Now let $<_{\mathcal{O}}$ be any ordering on the category $\mathcal{O}$ as in Subsection \ref{O}. 

\begin{prop}\label{can bas set}  Set ${\bf B}:= \{ E \in \Irr(W)\,|\, \KZ({\rm L}(E)) \neq 0\}.$
 \begin{enumerate}[(a)] 
 \item The set
$\{ \KZ({\rm L}(E)) \,|\,E \in {\bf B} \}$ is a complete set of pairwise non-isomorphic simple $\mathcal{H}_{{\bf k}}(W)$-modules. \smallbreak
\item For all $E \in \Irr(W),\,E' \in {\bf B}$, we have $[\Delta(E):{\rm L}(E')]=[\KZ(\Delta(E)):\KZ({\rm L}(E'))]$.\smallbreak
\item If $E\in {\bf B}$, then $[\KZ(\Delta(E)):\KZ({\rm L}(E))]=1$.\smallbreak
\item If $[\KZ(\Delta(E)):\KZ({\rm L}(E'))] \neq 0$ for some $E \in \Irr(W)$ and $E' \in {\bf B}$, then either 
$E=E'$ or $E' <_{\mathcal{O}} E$.
\end{enumerate}
\end{prop}

Property (a) follows from \cite[Theorem 5.14]{GGOR}. For the proof of  properties (b), (c) and (d), all of them deriving from the fact that $\KZ$ is exact, the reader may refer to \cite[Proposition 3.1]{CGG}.

The simple modules killed by the $\KZ$-functor are exactly the ones that do not have full support. Their determination, and thus the determination of the set ${\bf B}$, is a very difficult problem.

We also obtain a  decomposition matrix $D_{\bf k}$ for the specialised Hecke algebra $\mathcal{H}_{{\bf k}}(W)$ with respect to the specialisation  $\Theta$. The rows of $D_{\bf k}$ are indexed by $\Irr(W)$ and its columns by $\Irr(\mathcal{H}_{{\bf k}}(W))$. Following Proposition \ref{can bas set}, $D_{\bf k}$ can be obtained from the decomposition matrix ${\bf D}$ of the category $\mathcal{O}$ by removing the columns that correspond to the simple modules killed by the $\KZ$-functor, that is, the columns labelled by $\Irr(W) \setminus {\bf B}$. This implies that $D_{\bf k}$ becomes lower unitriangular when its rows are ordered with respect to $<_{\mathcal{O}}$, in the same way that, in the cases where $\Theta$ factors through a cyclotomic Hecke algebra, the existence of a canonical basic set implies that $D_{\bf k}$ becomes lower unitriangular when its rows are ordered with respect to the $a$-function. If we could show that the $a$-function defines an ordering on the category $\mathcal{O}$, we would automatically obtain the existence of a canonical basic set for  $\mathcal{H}_{{\bf k}}(W)$. At the same time, we would obtain the determination of ${\bf B}$ in the cases where canonical basic sets have  already been  explicitly described.

\subsection{The $(a+A)$-function}

Let  ${\bf m}=(m_{\mathcal{C},j})_{(\mathcal{C} \in \mathcal{A}/W)(0 \leq j \leq e_{\mathcal{C}}-1)}$  be a family of integers and let  ${\varphi_{\textbf{m}}}: v_{\mathcal{C},j} \mapsto q^{m_{\mathcal{C},j}}$ be the corresponding cyclotomic specialisation for the Hecke algebra $\mathcal{H}(W)$. Let  $\theta: q \mapsto \eta$ be a specialisation such that $\eta$ is a non-zero complex number. If $\eta$ is not a root of unity or $\eta=1$, then, due to Theorem \ref{sscriterion}  and the form of the Schur elements of $\mathcal{H}_{{\varphi_{\textbf{m}}}}(W)$, the specialised Hecke algebra $\mathcal{H}_\eta(W)$ is semisimple. So we may assume from now on that $\eta$ is a root of unity of order $e>1$, namely $\eta=\zeta_e^r$ for some $r \in \Z_{>0}$ such that ${\rm gcd}(e,r)=1$.

Let  ${\bf k}=(k_{\mathcal{C},j})_{(\mathcal{C} \in \mathcal{A}/W)(0 \leq j \leq e_{\mathcal{C}}-1)}$  be the family of rational numbers defined by 
$$
k_{\mathcal{C},j}:=\frac{rN_W}{e}\, m_{\mathcal{C},j} \,\,\,\,\,\text{for all }\,\mathcal{C},j.
$$
Then $\mathcal{H}_{\bf k}(W) = \mathcal{H}_\eta(W)$. Following \cite[\S 3.3]{CGG}, we obtain the following equation which relates the functions $a^{\bf m}$ and $A^{\bf m}$ for $\mathcal{H}_{{\varphi_{\textbf{m}}}}(W)$ with the $c$-function for ${\bf H}_{{\bf k}}(W)$:

 \begin{equation} \label{a+A}
 a^{\bf m}_E+A^{\bf m}_E= \frac{e}{rN_W} c_E + \sum_{H\in \mathcal{A}}\sum_{j=0}^{e_\mathcal{C}-1} m_{\mathcal{C},j}\,\,\,\,\, \text{ for all }\, E \in \Irr(W),
 \end{equation}
where $\mathcal{C}$ denotes the orbit of $H \in \mathcal{A}$ under the action of $W$. 

\begin{rem}{\rm The above formula was also obtained in \cite[\S 6.2]{GGOR} for finite Weyl groups in the equal parameter case.}
\end{rem}

Equation (\ref{a+A}) implies that $a^{\bf m}+A^{\bf m}$ yields the same ordering on $\Irr(W)$ as the $c$-function (note that in this case $c_E \in \Q$ for all $E \in \Irr(W)$). Thus, 
 $a^{\bf m}+A^{\bf m}$ is also an ordering on the category $\mathcal{O}$, that is,  if $[\Delta(E):{\rm L}(E')] \neq 0$ for some 
$E,E' \in \Irr(W)$, then either $E=E'$ or $a^{\bf m}_{E'}\,+\,A^{\bf m}_{E'}\,< \, a^{\bf m}_E\,+\,A^{\bf m}_E$. 
If now the function $a^{\bf m}$ is compatible with $a^{\bf m}+A^{\bf m}$, that is,  for all $E,E' \in \Irr(W)$, 
\begin{equation}\label{aA2}
 a^{\bf m}_{E'}\,+\,A^{\bf m}_{E'} \,<\,a^{\bf m}_E\,+\,A^{\bf m}_E \,\,\Rightarrow\,\, a^{\bf m}_{E'}\,<\,a^{\bf m}_{E},
 \end{equation}
 then $a^{\bf m}$ is an ordering on the category $\mathcal{O}$ and we obtain the existence of a canonical basic set for $\mathcal{H}_{{\varphi_{\textbf{m}}}}(W)$ with respect to $\theta$ by Proposition \ref{can bas set}. This is true in several cases, but unfortunately not true in general. Some exceptional complex reflection groups where (\ref{aA2}) holds and the above argument works are: 
 \begin{center}
 $G_{23}=H_3$,\, $G_{24}$,\, $G_{27}$,\, $G_{29}$ \,and \,$G_{30}=H_4$ \footnote{The groups $G_{23}$, $G_{24}$, $G_{27}$, $G_{29}$, $G_{30}$, $G_{31}$, $G_{33}$, $G_{34}$, $G_{35}$, $G_{36}$ and $G_{37}$ are easy to check with a computer; they are all generated by pseudo-reflections of order $2$ whose reflecting hyperplanes belong to the same orbit.}. 
 \end{center}
 This yields the existence of canonical basic sets for the groups $G_{24}$, $G_{27}$ and $G_{29}$, which was not known before. To summarise, we have the following:
 
 \begin{prop}
 Let $W=G_n$, $n \in \{23,24,27,29,30\}$. Let ${\bf m}$ and ${\bf k}$ be defined as above, and let $E,\,E' \in \Irr(W)$.
  If $[\Delta(E):{\rm L}(E')] \neq 0$, then either $E=E'$ or $a^{\bf m}_{E'}\,< \, a^{\bf m}_E$. In particular, we have
 $\KZ({\rm L}(E)) \neq 0$ if and only if $E$ belongs to the canonical basic set of $\mathcal{H}_{{\varphi_{\textbf{m}}}}(W)$ with respect to $\theta: q \mapsto \zeta_e^r$.
 \end{prop}
 
 \subsection{Canonical basic sets for Iwahori--Hecke algebras from rational Cherednik algebras}
 Equation (\ref{a+A}) has also allowed us to show that, in the case where $W$ is a finite Coxeter group, and  assuming that Lusztig's conjectures {\bf P1 -- P15} hold,  the $c$-function is compatible with the ordering $\leq_{\mathcal{L}\mathcal{R}}$ on two-sided cells, since $a$ and $A$ are (see \cite[Remark 5.4]{Geck2} for the $a$-function, \cite[Corollary 21.6]{31} and \cite[Proposition 2.8]{chja1} for $A$).  This in turn was  crucial in showing \cite[Corollary 4.7]{CGG}:

 \begin{prop}\label{propcanfromCh}
Let $(W,S)$ be a finite Coxeter group and let $\mathcal{H}(W,L)$ be the Iwahori--Hecke algebra of $W$ with parameter $L$, as defined in \S\ref{KLCells}.
For $H \in \mathcal{A}$, let $s_H \in W$ be the reflection with reflecting hyperplane $H$ and let $\mathcal{C}$ be the orbit of $H$ under the action of $W$. 
If $H' \in \mathcal{C}$, then we have $L(s_H)=L(s_{H'})$ and we can set $L_{\mathcal{C}}:=L(s_H)$.
Let $e, r \in \mathbb{Z}_{>0}$ such that ${\rm gcd}(e,r)=1$, and take, for all $\mathcal{C} \in \mathcal{A}/W$,
$$ k_{\mathcal{C},0} = \frac{r L_{\mathcal{C}}}{e} \,\,\,\, \text{ and }\,\,\,\,\, k_{\mathcal{C},1} = -\frac{rL_{\mathcal{C}}}{e} .$$
If $E \in \Irr(W)$, then $\KZ({\rm L}(E)) \neq 0$ if and only if $E$ belongs to the canonical basic set of $\mathcal{H}(W,L)$ with respect to $\theta: q \mapsto \zeta_e^r$.
 \end{prop}
 
 The proof uses a connection,  established in  \cite[Proposition 4.6]{CGG}, between category $\mathcal{O}$ and the cellular structure of the Iwahori--Hecke algebra. More specifically, if $E \in \Irr(W)$, then $\KZ(\Delta(E))$ is isomorphic to the cell module $W_\theta(E)$ defined in \cite[Example 4.4]{cellular};
we will not go into further details here. Note though that, in Proposition \ref{propcanfromCh}, we have not included the assumption that Lusztig's conjectures must hold. The reason is that the only case where they are not known to hold, the case of $B_n$,  is covered by Corollary \ref{gl1n} below.

 \begin{rem} {\rm The above result can be generalised to the case where $k_{\mathcal{C},0} = \lambda L_{\mathcal{C}}$ and
  $k_{\mathcal{C},1} = - \lambda L_{\mathcal{C}}$for any complex number $\lambda$. If $\lambda \in \Z$ or $\lambda \in \C \setminus \Q$, then both category $\mathcal{O}$ and $\mathcal{H}_{{\rm exp}(2\pi i \lambda)}(W,L)$ are semisimple, so the statement trivially holds.
 If $\lambda$ is a negative rational  number, let us say $\lambda = - r /e$ for some  $e, r \in \mathbb{Z}_{>0}$ with ${\rm gcd}(e,r)=1$, and $E \in \Irr(W)$, then $\KZ({\rm L}(E)) \neq 0$ if and only if $E$ belongs to the canonical basic set of $\mathcal{H}(W,-L)$ with respect to $\theta: q \mapsto \zeta_e^r$. We recall now that the canonical basic sets for finite Coxeter groups where $L$ can take negative values are described in \cite{chja1}. In fact, $E$ belongs to the canonical basic set of $\mathcal{H}(W,-L)$ with respect to $\theta: q \mapsto \zeta_e^r$ if and only if  $E \otimes \varepsilon$ belongs to the canonical basic set of $\mathcal{H}(W,L)$ with respect to $\theta$, where $\varepsilon$ denotes the sign representation of $W$. 
 } 
 \end{rem}
 
 Proposition \ref{propcanfromCh} yields the existence of canonical basic sets for all finite Coxeter groups in a uniform way. At the same time, it yields a description of the simple modules that are not killed by the $\KZ$-functor, since canonical basic sets for finite Coxeter groups are explicitly known (see, for example, \cite{GJlivre}).
However, it does not imply that the $a$-function is an ordering on the category $\mathcal{O}$, because we do not know what happens with the simple modules killed by the $\KZ$-functor.
We do believe though that, for finite Coxeter groups,  the $a$-function is an ordering on the category $\mathcal{O}$.
 
\begin{exmp}{\rm Let $W$ be the symmetric group $\mathfrak{S}_n$ and let $l:=L(s)$ for every transposition $s \in \mathfrak{S}_n$ (there exists only one orbit $\mathcal{C}$ in $\mathcal{A}/W$). 
Let $\eta^{2l}:=\zeta_e^{r}$ for some  $e, r \in \mathbb{Z}_{>0}$ with ${\rm gcd}(e,r)=1$. As we saw in Example \ref{excanbassym}, the canonical basic set $\mathcal{B}_\theta$  of $\mathcal{H}(W,l)$ with respect to $\theta: q \mapsto \eta$ consists of the $e$-regular partitions of $n$.
Now take $k_{\mathcal{C},0}=r/2e$ and  $k_{\mathcal{C},1}=-r/2e$.  Let $\lambda$ be a partition of $n$ and let $E^\lambda$ be the corresponding irreducible representation of $\mathfrak{S}_n$. We have  $\KZ({\rm L}(E^\lambda)) \neq 0$ if and only if $\lambda$ is $e$-regular.
}
\end{exmp}

\subsection{Canonical basic sets for Ariki--Koike algebras from rational Cherednik algebras}
 As we have said and seen earlier, there exist several orderings on the category $\mathcal{O}$. For $W=G(\ell,1,n)$, where the irreducible representations are parametrised by the $\ell$-partitions of $n$, one combinatorial ordering on the category $\mathcal{O}$ is given by Dunkl and Griffeth in \cite[Theorem 4.1]{DG}. More precisely, 
 in this case, 
 there are two hyperplane orbits in $\mathcal{A}/W$; we will denote them by $\mathcal{C}_ \textbf{s}$ and $\mathcal{C}_ \textbf{t}$. We have $e_{\mathcal{C}_ \textbf{s}}= \ell$ and $e_{\mathcal{C}_ \textbf{t}}= 2$. Let $(s_0, \ldots , s_{\ell - 1}) \in \Z^{\ell}$ and $e \in \Z_{>0}$.
 We define
${\bf k}=(k_{\mathcal{C}_ \textbf{s},0},\ldots,k_{\mathcal{C}_ \textbf{s},\ell-1},k_{\mathcal{C}_ \textbf{t},0}, k_{\mathcal{C}_ \textbf{t},1})$ by
\begin{equation}\label{params}
k_{\mathcal{C}_ \textbf{s},j}= \frac{s_j}{e} - \frac{j}{\ell}\,\,\,\, \text{ for }\, j =0,\ldots,\ell-1 ,\,\,\,\, k_{\mathcal{C}_ \textbf{t},0}=\frac{1}{e},\,\,\,k_{\mathcal{C}_ \textbf{t},1}=0.
\end{equation}
Then the $\KZ$-functor goes from the category $\mathcal{O}$ for ${\bf H}_{\bf k}(W)$ to the category of representations of
  the specialised Ariki--Koike algebra  $\mathcal{H}_{\bf k}(W)$ with relations 
$$
(\textbf{s} - \zeta_e^{s_0})(\textbf{s} - \zeta_e^{s_1})\cdots (\textbf{s}- \zeta_e^{s_{\ell - 1}}) = 0, \qquad 
(\textbf{t}_i-\zeta_e)(\textbf{t}_i+1) = 0\,\,\,\,\text{ for }i=1,\ldots,n-1,
$$
 as in (\ref{ArikiKoike}).
 
Let $\bl =(\lambda^{(0)},\ldots,\lambda^{(\ell-1)})$ be an $\ell$-partition of $n$.  We will denote by $E^{\bl}$  the corresponding irreducible representation of  $G(\ell,1,n)$.  We define the set of nodes of $\bl$ to be  the set
$$[\bl]=\{ (a,b,c): 0 \leq c \leq \ell-1,\,\,a \geq 1,\,\,1 \leq b \leq \lambda_a^{(c)}\}.$$
Let $\gamma=(a(\gamma),b(\gamma),c(\gamma)) \in [\bl]$. We set $\vartheta(\gamma):=b(\gamma)-a(\gamma)+s_{c(\gamma)}$.
We then have the following  \cite[Proof of Theorem 4.1]{DG}:

\begin{prop}\label{DuGr} Let $\bl,\,\bl'$ be $\ell$-partitions of $n$.
 If
$[\Delta(E^\bl):{\rm L}(E^{\bl'})] \neq 0$, then there exist orderings $\gamma_1,\gamma_2,\ldots,\gamma_n$ and $\gamma_1',\gamma_2',\ldots,\gamma_n'$ of the nodes of $\bl$ and $\bl'$ respectively, and non-negative integers $\mu_1,\mu_2,\ldots,\mu_n$,  such that, for all $1\leq i\leq n$,
$$\mu_i \equiv c(\gamma_i)-c(\gamma_i')\,\,\mathrm{mod}\,\ell
\,\,\,\textrm{ and }\,\,\,
\mu_i= c(\gamma_i) - c(\gamma_i') + \frac{\ell }{e}(\vartheta(\gamma_i')-\vartheta(\gamma_i)).$$
\end{prop}

Now, there are several different cyclotomic Ariki--Koike algebras that produce the specialised Ariki--Koike algebra $\mathcal{H}_{\bf k}(W)$ defined above and they may have distinct $a$-functions attached to them. Using the combinatorial description of the $a$-function for $G(\ell,1,n)$ given in \cite[\S 5.5]{GJlivre}\footnote{This definition captures all $a$-functions  for $G(\ell , 1,n)$ in the literature: the function $a^{\bf m}$ for $m_{\mathcal{C}_ \textbf{s},j}=s_j\ell -ej$, $j=0,\ldots,\ell-1$,  given by Jacon \cite{jaca} and studied in the context of Uglov's work on canonical bases for higher level Fock spaces, and also the $a$-function for type $B_n$ ($\ell=2$) arising from the Kazhdan--Lusztig theory for  Iwahori--Hecke algebras with unequal parameters (see \cite[6.7]{GJlivre}).}, we showed in \cite[\S 5]{CGG} that  it is compatible  with the ordering on category $\mathcal{O}$ given by Proposition \ref{DuGr}. Consequently, the $a$-function also defines a highest weight structure on $\mathcal{O}$, that is, we have the following:

\begin{prop}
Let $\bl,\,\bl'$ be $\ell$-partitions of $n$. If
$[\Delta(E^\bl):{\rm L}(E^{\bl'})] \neq 0$, then either
 $\bl=\bl'$ or $a_{E^{\bl'}}<a_{E^\bl}$.
\end{prop}

The above result, combined with Proposition \ref{can bas set}, yields the following:

\begin{cor}\label{gl1n}
Let $W=G(\ell,1,n)$. Let $(s_0, \ldots , s_{\ell - 1}) \in \Z^{\ell}$ and $e \in \Z_{>0}$.
Let
${\bf k}=(k_{\mathcal{C}_ \textbf{s},0},\ldots,k_{\mathcal{C}_ \textbf{s},\ell-1},k_{\mathcal{C}_ \textbf{t},0}, k_{\mathcal{C}_ \textbf{t},1})$ be defined as in $(\ref{params})$.
If $\bl$ is an $\ell$-partition of $n$, then $\KZ({\rm L}(E^{\bl})) \neq 0$ if and only if $E^{\bl}$ belongs to the canonical basic set for $\mathcal{H}_{\bf k}(W)$ with respect to the $a$-function above.
\end{cor}

Thus, we obtain the existence of  canonical basic sets for Ariki--Koike algebras without the use of Ariki's Theorem.
On the other hand, the description of the canonical basic sets for Ariki--Koike algebras by \cite[Main Theorem]{jaca} yields a description of the set 
${\bf B}=\{ E^{\bl} \in \Irr(W)\,|\, \KZ({\rm L}(E^{\bl})) \neq 0\}$: we have that $E^{\bl} \in {\bf B}$ if and only if $\bl$ is an Uglov $\ell$-partition.

Finally, we expect a result similar to Corollary \ref{gl1n} to hold in the case where $W=G(\ell,p,n)$ for $p>1$.

\section{Rational Cherednik Algebras at $t = 0$}\label{lastsection}

Let us now consider the rational Cherednik algebra ${\bf H}_{0,\bc}(W)$. In this case, the centre of ${\bf H}_{0,\bc}(W)$ is isomorphic to the  spherical subalgebra of ${\bf H}_{0,\bc}(W)$, that is, $Z({\bf H}_{0,\bc}(W)) \cong e {\bf H}_{0,\bc} e$, where $e:=\frac{1}{|W|} \sum_{w\in W}w$.  So ${\bf H}_{0,\bc}(W)$ is a finitely generated $Z({\bf H}_{0,\bc}(W))$-module. From now on, we set $Z_\bc(W):=Z({\bf H}_{0,\bc}(W))$.

\subsection{Restricted rational Cherednik algebras}In the case of finite Coxeter groups the following was proved in \cite[Proposition 4.15]{EG}, and the general case is due to \cite[Proposition 3.6]{go}.

\begin{prop}\begin{enumerate}[(i)]
\item The subalgebra $\mathfrak{m}:=\C[\mathfrak{h}]^W \otimes \C[\mathfrak{h}^*]^W$  of ${\bf H}_{0,\bc}(W)$ is contained in $Z_\bc(W)$. \smallbreak
\item $Z_\bc(W)$ is a free $\mathfrak{m}$-module of rank $|W|$.
\end{enumerate}
\end{prop}

Let $\mathfrak{m}_+$ denote the ideal of $\mathfrak{m}$ consisting of elements with zero constant term.

\begin{defn}
{\rm We define the \emph{restricted rational Cherednik algebra} to be}
$$\overline{{\bf H}}_{0,\bc}(W):={\bf H}_{0,\bc}(W) /\mathfrak{m}_+ {\bf H}_{0,\bc}(W).$$
\end{defn}

This algebra was originally introduced, and extensively studied, in \cite{go}. The PBW Theorem implies that, as a vector space,
$$\overline{{\bf H}}_{0,\bc}(W) \cong \C[\mathfrak{h}]^{coW} \otimes \C W \otimes \C[\mathfrak{h}^*]^{coW}$$
where $\C[\mathfrak{h}]^{coW}=\C[\mathfrak{h}]/\langle \C[\mathfrak{h}]^W_+ \rangle$ 
is the \emph{coinvariant algebra}. Since $W$ is a complex reflection group, $\C[\mathfrak{h}]^{coW}$ has dimension $|W|$ and is isomorphic to the regular representation as a $\C W$-module. Thus, ${\rm dim}_{\C} \overline{{\bf H}}_{0,\bc}(W) = |W |^3$.

Let $E \in  \Irr(W)$. We set
$$\overline{\Delta}(E):=\overline{{\bf H}}_{0,\bc}(W)  \otimes_{\C[\mathfrak{h}^*]^{coW} \rtimes W}  E ,$$
where $\C[\mathfrak{h}^*]^{coW}$ acts trivially on  $E$ (that is, $\C[\mathfrak{h}^*]^{coW}_+$ acts on $E$ as zero)
and $W$ acts naturally. The module $\overline{\Delta}(E)$ is the \emph{baby Verma module} of $\overline{{\bf H}}_{0,\bc}(W)$ associated to $E$. We summarise, as is done in \cite[Proposition 4.3]{go}, the results of \cite{HoNa} applied to this situation. 

\begin{prop} Let $E,\,E' \in  \Irr(W)$.
\begin{enumerate}[(i)]
\item The baby Verma module $\overline{\Delta}(E)$ has a simple head, $\overline{{\rm L}}(E)$. Hence, $\overline{\Delta}(E)$ is indecomposable. \smallbreak
\item $\overline{\Delta}(E) \cong \overline{\Delta}(E')$ if and only if $E \cong E'$. \smallbreak
\item The set $\{\overline{{\rm L}}(E) \,|\, E \in  \Irr(W)\}$  is a complete set of pairwise non-isomorphic simple $\overline{{\bf H}}_{0,\bc}$-modules.
\end{enumerate}
\end{prop}

\subsection{The Calogero--Moser partition}
Recall that the generalised Calogero--Moser space $X_\bc(W)$  is defined to be the affine variety ${\rm Spec}\,Z_\bc(W)$. By Theorem \ref{smoothness}, $(\mathfrak{h} \oplus \mathfrak{h}^*)/W$ admits a symplectic resolution if and only if  $X_\bc(W)$ is smooth for generic values of $\bc$. Etingof and Ginzburg proved that  $X_{\bc}(G)$ is smooth for generic $\bc$ when $W=G(\ell,1,n)$
 \cite[Corollary 1.14]{EG}. Later, Gordon showed that $X_{\bc}(G)$ is a singular variety for all choices of the parameter $\bc$ for the following finite Coxeter groups \cite[Proposition 7.3]{go}: $D_{2n}\, (n \geq 2)$, $E_6$, $E_7$, $E_8$, $F_4$, $H_3$, $H_4$ and $I_2(m) \,(m \geq 5)$.

Now, since the algebra  $\overline{{\bf H}}_{0,\bc}$ is finite dimensional, we can define its blocks in the usual way (see \S \ref{definitionofblocks}).
Let $E,\, E'  \in  \Irr(W)$.
Following \cite{GM}, we define
the \emph{Calogero--Moser partition} of $\Irr(W)$ to be the set of equivalence classes of $\Irr(W)$ under the equivalence relation: 
\begin{center}
$E \sim_{CM} E'$ \,if and only if\, $\overline{{\rm L}}(E)$ and $\overline{{\rm L}}(E')$ belong to the same block.
\end{center}
We will simply write $CM_{\bc}$-partition for the Calogero--Moser partition of $\Irr(W)$. 
The inclusion $\mathfrak{m} \subset Z_\bc(W)$ defines a finite surjective morphism 
$$\mathcal{Y}: X_\bc(W) \longrightarrow \mathfrak{h}/W \times \mathfrak{h}^*/W$$
where $\mathfrak{h}/W \times \mathfrak{h}^*/W = {\rm Spec}\,\mathfrak{m}$. M\"uller's theorem (see \cite[Corollary 2.7]{BrGo}) implies that the natural map
$\Irr(W) \rightarrow \mathcal{Y}^{-1}(0),\, E \mapsto {\rm Supp}(\overline{{\rm L}}(E))$ factors through the $CM_{\bc}$-partition. Using this fact, one can show that the geometry of $X_\bc(W)$ is related to the $CM_{\bc}$-partition
in the following way.

\begin{thm}\label{parabolic} The following are equivalent:
\begin{enumerate}[(1)]
\item The generalised Calogero--Moser space $X_\bc(W)$ is smooth. \smallbreak
\item The  $CM_{\bc}$-partition of $\Irr(W')$ is trivial for every parabolic subgroup $W'$ of $W$.
\end{enumerate}
\end{thm}

Using the above result and the classification of irreducible complex reflection groups (see Theorem \ref{ShToClas}), Bellamy has shown the following  \cite[Theorem 1.1]{bel}:

\begin{thm} Let $W$ be an irreducible complex reflection group. The  generalised Calogero--Moser space $X_\bc(W)$ is smooth for generic values of $\bc$ if and only if $W$ is of type $G(\ell,1,n)$ or $G_4$. In every other case, $X_\bc(W)$ is singular for all choices of $\bc$.
\end{thm}

\begin{cor}
Let $W$ be an irreducible complex reflection group. The space  $(\mathfrak{h} \oplus \mathfrak{h}^*)/W$ admits a symplectic resolution if and only if  
$W$ is of type $G(\ell,1,n)$ or $G_4$. 
\end{cor}

\subsection{The Calogero--Moser partition and Rouquier families}
It just so happens that  the cases where $X_\bc(W)$ is generically smooth, and the Calogero--Moser partition generically trivial, are exactly the cases where the Rouquier families are generically trivial (that is, the Rouquier families associated with no essential hyperplane are singletons). This, combined with the fact that the Calogero--Moser partition into blocks enjoys some property of semicontinuity, led to the question whether there is a connection between the two partitions.

The question was first asked  by Gordon and Martino \cite{GM} in terms of a connection between the Calogero--Moser partition and families of characters for type $B_n$. In their paper, they computed the  $CM_{\bc}$-partition, for all $\bc$, for complex reflection groups of type $G(\ell,1,n)$ and showed that for $\ell=2$, using the conjectural combinatorial description of Kazhdan--Lusztig cells for type $B_n$ by \cite{BGIL},  the $CM_{\bc}$-partition coincides with the partition into Kazhdan--Lusztig families. After that, Martino  \cite{Mau} compared 
the combinatorial description of the  $CM_{\bc}$-partition for type $G(\ell,1,n)$ given in \cite{GM} with the description of the partition into Rouquier families, given by \cite{berkeley}, for a suitable cyclotomic Hecke algebra $\mathcal{H}_{\bc}$ of $G(\ell,1,n)$ (different from the one defined in \S\ref{KZ}). He showed that the two partitions coincide when $\ell$ is a power of a prime number (which includes the cases of type $A_n$ and $B_n$), but not in general. In fact, he showed that the 
 $CM_{\bc}$-partition for $G(\ell,1,n)$ is the same as the one obtained by \cite{BK}. He thus obtained the following two connections  between the $CM_{\bc}$-partition and the partition into Rouquier families for $G(\ell,1,n)$, and he conjectured that they hold for every complex reflection group $W$
\cite[2.7]{Mau}:
\begin{enumerate}[(a)]
\item The  $CM_{\bc}$-partition for generic $\bc$ coincides with the generic partition into Rouquier families (both being trivial for $W = G(\ell,1,n)$); \smallbreak
\item The partition into Rouquier families refines the $CM_{\bc}$-partition, for all choices of $\bc$; that is, if $E,E' \in \Irr(W)$ belong to the same Rouquier family of $\mathcal{H}_{\bc}$, then $E \sim_{CM} E'$.
\end{enumerate}

Conditions (a) and (b) are known as ``Martino's Conjecture''. Using the combinatorics of \cite{GM} and \cite{Mau},
Bellamy computed the  $CM_{\bc}$-partition, for all $\bc$, and proved Martino's conjecture in the case where $W$ is of type $G(\ell,p,n)$  \cite{bellpn}; note that when $p>1$ the generic partitions in this case are not trivial.  However, a counter-example for (a)  was found recently by Thiel \cite{Thiel} in the case where $W=G_{25}$.
Thiel calculated the $CM_{\bc}$-partition for generic $\bc$  for the exceptional complex reflection groups $G_4,\, G_5,\,G_6,\,G_8,\, G_{10}, G_{23}=H_3, G_{24}, G_{25}$ and $G_{26}$. Comparing his results with the generic partition into Rouquier families for these groups, given by \cite{mybook}, he showed that Part (a) of Martino's Conjecture holds in every case\footnote{for $G_4$ this was already known by \cite{bel}.} except for when $W=G_{25}$. In this particular case, the generic partition into Rouquier families simply refines the $CM_{\bc}$-partition for generic $\bc$. So we will state here as a conjecture only Part (b) of Martino's conjecture, which is still an open problem, and proved in all the above cases.

\begin{conj}{\rm  \textbf{(Martino's Conjecture)} } Let $W$ be a complex reflection group. The partition into Rouquier families (for a suitably chosen cyclotomic Hecke algebra $\mathcal{H}_{\bc}$ of $W$) refines the $CM_{\bc}$-partition, for all choices of $\bc$; that is, if $E,E' \in \Irr(W)$ belong to the same Rouquier family of $\mathcal{H}_{\bc}$, then $E \sim_{CM} E'$.
\end{conj}

\begin{rem}{\rm
Note that, in all the cases checked so far where $W$ is a finite Coxeter group, the partition into Rouquier families and the $CM_{\bc}$-partition coincide. This covers the finite Coxeter groups of types $A_n$, $B_n$, $D_n$ and the dihedral groups for all choices of $\bc$, and $H_3$ for generic $\bc$.
}
\end{rem}

\subsection{The Calogero--Moser partition and Kazhdan--Lusztig cells}
 In an effort to develop a  generalised Kazhdan--Lusztig cell theory, Bonnaf\'e and Rouquier used the Calogero--Moser partition to define, what they call in \cite{Bora}, \emph{Calogero--Moser cells} for all complex reflection groups. An advantage of this, quite geometric, approach is that the Calogero--Moser partition exists naturally for all complex reflection groups.
It also implies automatically the existence of a semicontinuity property for cells, a property that was conjectured and proved in some cases for Kazhdan--Lusztig cells by Bonnaf\'e \cite{semicont}. However, Calogero--Moser cells are very hard to compute and their construction depends on an ``uncontrollable'' choice. After very long computations by Bonnaf\'e and Rouquier,  it is now confirmed that the Calogero--Moser cells coincide with the Kazhdan--Lusztig cells in the smallest possible cases ($A_2$, $B_2$, $G_2$);  there is still a lot of work that needs to be done.

\section{Acknowledgements}
 First of all, I would like to thank the Mathematical Sciences Research Institute (MSRI) for its hospitality during the programme ``Noncommutative Algebraic Geometry and Representation Theory'' and for its support.  My thanks to the organisers of the Introductory Workshop, Michael Artin, Michel Van den Bergh and Toby Stafford, who invited me to give two talks on symplectic reflection algebras; the current paper is inspired by these talks. A special thanks to Toby Stafford for giving me more time to finish it. I feel a lot of gratitude towards Gwyn Bellamy and Iain Gordon for answering all my questions on symplectic reflection algebras (and they were a lot!), and for writing two excellent surveys on the topic, \cite{MSRI} and \cite{ICRA}. Iain also offered to read this manuscript and suggested many corrections, for which I am grateful. I also thank C\'edric Bonnaf\'e for answering my questions in Luminy, reading this manuscript and making useful comments.
My deepest gratitude towards Guillaume Pouchin for going through all the phases of this project with me. 
Finally, I would like to thank Gunter Malle and the referee of this paper for suggesting many corrections to the final version.

This material is based upon work supported by the National Science Foundation under Grant No. 0932078 000, while the author was in residence at the Mathematical Science Research Institute (MSRI) in Berkeley, California, during 2013.
This research project is also implemented within the framework of the Action ``Supporting Postdoctoral Researchers'' of the Operational Program
``Education and Lifelong Learning'' (Action's Beneficiary: General Secretariat for Research and Technology), and is co-financed by the European Social Fund (ESF) and the Greek State.

\end{document}